\documentclass[preprint,12pt]{elsarticle}

\usepackage{a4wide,latexsym,varioref,amsmath,bbm,amsfonts,xcolor}

\newcommand{\numsection}[1]{\section{#1}\setcounter{equation}{0}}

\renewcommand{\theequation}{\arabic{section}.\arabic{equation}}

\newcommand{\calA}{{\cal A}}
\newcommand{\calM}{{\cal M}}
\newcommand{\calN}{{\cal N}}
\newcommand{\calS}{{\cal S}}
\newcommand{\calX}{{\cal X}}
\newcommand{\beqn}[1]{\begin{equation}\label{#1}}
\newcommand{\eeqn}{\end{equation}}
\newcommand{\ii}[1]{\{1, \ldots, #1 \}}
\newcommand{\iibe}[2]{\{ #1, \ldots, #2 \}}
\newcommand{\flow}{f_{\rm low}}
\newcommand{\barphi}{\overline{\phi}}
\newcommand{\DT}{\Delta T}
\newcommand{\barDT}{\overline{\Delta T}}
\newcommand{\algn}{{\footnotesize {\sf IAR$qp$} }}
\newcommand{\expect}{\mathbb{E}}
\newcommand{\indic}{\mathbbm{1}_}
\newcommand{\no}[1]{#1^c}
\newcommand{\AM}{\calM}
\newcommand{\indicM}{\indic{\AM_k}}
\newcommand{\indicMc}{\indic{\no{\AM_k}}}
\newcommand{\indicS}{\indic{\calS_k}}
\newcommand{\indicSc}{\indic{\no{\calS_k}}}
\newcommand{\prob}{\mathbb{P}{\rm r}}
\newcommand{\sigax}{\calA_{k-1}^{M}}
\newcommand{\pM}{p_*}
\newcommand{\sfrac}[2]{{\scriptstyle \frac{#1}{#2}}}
\newcommand{\half}{\sfrac{1}{2}}
\newcommand{\tim}[1]{\;\; \mbox{#1} \;\;}
\newcommand{\req}[1]{(\ref{#1})}
\newcommand{\eqdef}{\stackrel{\rm def}{=}}
\newcommand{\ms}{\;\;\;\;}
\newcommand{\bard}{\overline{d}}
\newcommand{\barf}{\overline{f}}
\newcommand{\barD}{\overline{D}}

\newcommand{\barT}{\overline{T}}
\newcommand{\bigfrac}[2]{\frac{\displaystyle #1}{\displaystyle #2}}
\newcommand{\bigsum}{\displaystyle \sum}
\newcommand{\bigmax}{\displaystyle \max}
\newtheorem{theorem}{Theorem}[section]
\newtheorem{lemma}[theorem]{Lemma}

\newtheorem{definition}{Definition}
\newcommand{\llem}[2]{\vspace{\baselineskip} 
\noindent\framebox[\textwidth]{\parbox{0.95\textwidth}{
\begin{lemma} \label{#1} \rm #2 \end{lemma} } } \vspace{\baselineskip} }
\newcommand{\lthm}[2]{\vspace{\baselineskip} 
\noindent\framebox[\textwidth]{\parbox{0.95\textwidth}{
\begin{theorem} \label{#1} \rm #2 \end{theorem} } } \vspace{\baselineskip} }
\newcommand{\bpr}{{\bf Proof.} \hspace{1.5mm}}
\newcommand{\epr}{\hfill $\Box$ \vspace*{1em}}
\newcommand{\proof}[1]{
\begin{list}{}{
\setlength{\topsep}{0.0pt}
\setlength{\partopsep}{0.0pt}
\setlength{\leftmargin}{0.025\textwidth}
\setlength{\rightmargin}{0.5\leftmargin}
\setlength{\labelwidth}{0.5\leftmargin}
\setlength{\labelsep}{0.25\leftmargin}}
\item \bpr #1 \epr \noindent
\end{list}}
\newcounter{algo}[section]
\renewcommand{\thealgo}{\thesection.\arabic{algo}}
\newcommand{\algo}[3]{\refstepcounter{algo}
\begin{center}\begin{figure}[htbp]
\framebox[\textwidth]{
\parbox{0.95\textwidth} {\vspace{\topsep}
{\bf Algorithm \thealgo : #2}\label{#1}\\
\vspace*{-\topsep} \mbox{ }\\
{#3} \vspace{\topsep} }}
\end{figure}\end{center}}
\renewcommand{\Re}{\hbox{I\hskip -2pt R}}
\newcommand{\smallRe}{\hbox{\footnotesize I\hskip -2pt R}}
\DeclareMathOperator*{\argmax}{arg\,max}

\begin{document}
\begin{frontmatter}
\title{Adaptive Regularization for Nonconvex Optimization
  Using Inexact Function Values and Randomly Perturbed Derivatives}

\author[SB]{ S. Bellavia}
\author[GB] {G. Gurioli}
\author[BM] {B. Morini}
\author[PT] {Ph. L. Toint} 
\address[SB] {Dipartimento di Ingegneria Industriale,
    Universit\`{a} degli Studi di Firenze, Italy. Member of the INdAM Research
    Group GNCS. Email: stefania.bellavia@unifi.it}
\address[GB]{Dipartimento di Matematica e Informatica ``Ulisse Dini'',
    Universit\`{a} degli Studi di Firenze, Italy.  {Member of the INdAM Research
    Group GNCS.} Email: gianmarco.gurioli@unifi.it}
\address[BM]{Dipartimento di Ingegneria Industriale,
    Universit\`{a} degli Studi di Firenze, Italy. Member of the INdAM Research
    Group GNCS. Email: benedetta.morini@unifi.it} 
\address[PT] { Namur Center for Complex Systems (naXys),
    University of Namur, 61, rue de Bruxelles, B-5000 Namur, Belgium.
    Email: philippe.toint@unamur.be}


\begin{abstract}
  A regularization algorithm allowing random noise in
  derivatives and inexact function values is proposed for computing
  approximate local critical points of any order for smooth unconstrained
  optimization problems. For an objective function with Lipschitz continuous
  $p$-th derivative and given an arbitrary optimality order $q \leq p$,
  it is shown that this algorithm will, in expectation, compute such a point in at most 
  $O\Bigl(\left(\min_{j\in\ii{q}}\epsilon_j\right)^{-\frac{p+1}{p-q+1}}\Bigr)$
  inexact evaluations of $f$ and its derivatives whenever $q\in\{1,2\}$, where
  $\epsilon_j$ is the tolerance for $j$th order accuracy.  This bound becomes
  at most $O\Bigl(\left(\min_{j\in\ii{q}}\epsilon_j\right)^{-\frac{q(p+1)}{p}}\Bigr)$
  inexact evaluations if $q>2$ and all derivatives are Lipschitz continuous.
  Moreover these bounds are sharp in the order of the accuracy tolerances.
  An extension to convexly constrained problems is also outlined.
\end{abstract}

 \end{frontmatter}
 
\numsection{Introduction}

We consider the evaluation complexity of an adaptive regularization algorithm for
computing approximate \emph{local} minimizers of arbitrary order for the
unconstrained minimization problem of the form
\beqn{problem}
\min_{x \in \smallRe^n}  f(x), 
\eeqn
where the objective function $f$ is sufficiently smooth and the values of its
$j$-th derivatives $\nabla_x^j f$ are subject to random noise and can only be computed inexactly.  
Inexact values of the objective function are also allowed, but their
inaccuracy is assumed to be deterministically controllable.

\vskip 5pt 
\noindent
\textbf{Motivation and context.} Without further comments, this statement of
the paper's purpose may be difficult to interpret, and we start by clarifying
the notion of ``evaluation complexity''.  Since the computational cost of
nonlinear optimization algorithms for the local solution of \req{problem} is
typically dominated by that of evaluating $f(x)$ and its derivatives, and
since these evaluations are fully independent of the algorithms (and unknown
to the algorithm designers), classical concepts of total computational
complexity are difficult to apply to such algorithms, except for the simplest
cases\footnote{For example, consider the typical
  optimization problem in deep learning applications, where the cost of the
  objective function depends on the size, depth and structure of the neural
  network used, all things that are not known by the algorithm.  Or
  optimization of large inverse problems like weather forecasting, where
  evaluating the objective function depends on solving a complicated
  multidimensional time-dependent partial-differential equation, whose dimension, domain
  shape, horizon, level of discretization and nonlinearity are also unknown and may
  vary from application of the algorithm to the next.}. This difficulty is at the origin of the now standard concept of
``evaluation complexity'' (sometimes called ``oracle complexity'') where the
cost of running an algorithm is approximated by the total cost of evaluating
the objective function and its relevant derivatives (the oracle), which is
then measured by counting the total number of such evaluations. This concept
has a long history in the optimization research community, and has generated a
vast literature covering many kinds of algorithms and problems types (see
\cite{NemiYudi83, Vava92, Vava93, Nest04, NestPoly06, CartGoulToin10a,
  CartGoulToin11d, BianChenYe17, CurtRobi18} and the reference therein for a
limited sample of this vibrant research area). Of course, like the classical
notion of complexity, this approach only applies when a class of problems is
well identified and when the algorithm's termination rules are clear. In our
case, these rules need to reflect the notion of approximate local solution of
problem \req{problem}, and will be detailed in due course (in
Section~\ref{algo-s}).

The class of problems of interest here (as summarized above) is the
approximate local minimization of smooth functions whose values (and that of
its derivatives) can only be computed inexactly.  Before providing more
technical detail on exactly what this means, it is useful to briefly review
the existing contributions in this domain\footnote{We focus here on algorithms
whose definition only involves quantities that are known to the user. In
particular, we explictly avoid requiring the knowledge of the problem
Lipschitz constants for the definition of the algorithm.}. Indeed, solving
optimization problems involving inexact evaluations is not a new topic and has
already been investigated in two different frameworks.  The first is that of
(deterministic) explicit dynamic accuracy, where it is assumed that the
accuracy of $f$ and its derivatives can be controlled by the algorithm (see
\cite[Section~10.6]{ConnGoulToin00}, \cite{KourHeinRizdvanB14},
\cite{GratSimoToin20} or \cite{BellGuriMoriToin19} for example).  In this
context, accuracy requirement are proposed that guarantee convergence to
approximate local solutions, and evaluation complexity of the resulting
algorithms can be analyzed \cite{BellGuriMoriToin19}, indicating a very modest
degradation of the worst-case performance compared with the case where
evaluations are exact \cite{CartGoulToin20b,CartGoulToin20a}.  A drawback of
this approach is that nothing is said for the case where the requested
accuracy requirement cannot be met or, as is often the case, cannot even be
measured. This problem does not occur in the second framework, in which the
inexactness in the function (and possible derivatives) values can be seen as
caused by some random noise, in which case the algorithm/user is not able to
specify an accuracy level and poor derivative approximations might result. The
available analysis for this case differ by the assumptions made on the
distribution of this noise. In \cite{PaquSche20}, the authors consider the
unbiased case and estimate the evaluation complexity a linesearch method for
finding approximate first-order critical points. Here, we assume  that
derivatives values can be approximated within a prescribed accuracy with a
fixed, sufficiently high probability, conditioned to the past.  A similar
context is considered in \cite{BeraCaoSche19}, where the objective function
values are inexact but computed with accuracy guaranteed with probability
one, and in \cite{PaquSche20}, where the authors consider the
unbiased case. A trust-region method (see \cite{ConnGoulToin00} for a full coverage of
such methods) is also proposed in \cite{ChenMeniSche18}, where it is proved to
converge almost-surely to first-order critical points.  Using similar
assumptions, the approach of \cite{CartSche17} includes the use of random
first-order models and directions within line search method as well as
probabilistic second-order models in Adaptive Cubic Regularization (ARC)
algorithms.  In both cases, the authors employ exact function evaluations. A
general theory for global convergence rate analysis is also provided. More
recently, \cite{BlanCartMeniSche19} proposed an evaluation complexity analysis
for a trust-region method (covering convergence to second-order points) using
elegant properties of sub-martingales, and making no assumption on bias. A
recent overview of this active research area is proposed in \cite{CurtSche20}.

\vskip 5pt 
\noindent
\textbf{Contributions.}  As suggested above, this paper deliberately considers
  the well-established concept of evaluation complexity, measuring the number
  of calls to user-supplied procedures for computing approximate function and
  derivatives values, irrespective of internal computations within the
  algorithm itself.  Of course,
  the authors are fully aware that total computational complexity (as opposed
  to evaluation complexity) is a different question. Fortunately, the
  difference is well-understood when searching for
  first- or second-order approximate local minimizers, in that moderately costly
  methods are available for handling the algorithm's internal calculations
  (see \cite{CarmDuch19, CartGoulToin12b, RoyeONeiWrig20}).  In other
  situations, the results presented here give an admittedly idealized but
  hopefully interesting estimation.

Having set the scene and clarified our objective, we now make the
contributions of this paper more precise.

\begin{itemize}
\item
  We consider finding approximate local solution of problem \req{problem} assuming that the objective function
  values can be computed within a prescribed accuracy, while at the same time allowing randomly inexact
  evaluations of its derivatives, thereby using a mix of the two frameworks
  described above. This work thus extends the analysis  provided in \cite{BellGuriMori21, BellGuriMoriToin19, XuRoosMaho20} 
  for adaptive regularization algorithms.

  Our assumptions on the type of inaccuracy allowed for the objective function
  complements that of \cite{BellGuri20, CartSche17}, allowing for more inexactness, but
  in a deterministic context.  This is a  realistic request in applications
  such as those where the objective function value is approximated by using
  smoothing operators and the derivatives are approximated by randomized
  finite differences \cite{BeraCaoSche19, BeraCaoChoSche20,MaggWachDoliStau18,NestSpok17}.
\item
  As in \cite{BellGuriMoriToin19,CartGoulToin20b,CartGoulToin20a}, we propose
  a regularization algorithm which is based on polynomial models of arbitrary
  degree. This obviously allows us to seek for first- and second-order
  critical points, as is standard, but we may also seek critical points of
  arbitrary order (we define what we mean by that in Section~\ref{algo-s}).
  In this respect we improve upon the algorithms with stochastic models such
  as \cite{CartSche17,BellGuri20, BellGuriMoriToin20, BlanCartMeniSche19, CurtSche20}.
\item
  We establish sharp worst-case bounds (in expectation) on the evaluation
  complexity of computing these (possibly high-order) approximate critical
  points, depending on the order and on the degree of the polynomial model
  used. Remarkably, these bounds correspond in order to the best known bounds for
  regularization algorithms using exact evaluations.
\end{itemize}
These results are obtained by a novel combination of the probabilistic framework of
\cite{CartSche17}, the approximation results of \cite{BellGuriMoriToin19} and
the proof techniques of \cite{CartGoulToin20a}.

\vskip 5pt 
\noindent
{\bf Outline.} The paper is organized as follows.  Section~\ref{algo-s} discusses
optimality measures for arbitrary order and introduces the regularization
algorithm and the associated probabilistic assumptions. Its evaluation complexity
is then studied in Section~\ref{complexity-s}.  
We finally present some conclusions and perspectives in Section~\ref{concl-s}.

\vskip 5pt 
\noindent
{\bf Notations.} Unless otherwise specified, $\|\cdot\|$ denotes the standard
Euclidean norm for both vectors and matrices. For a general symmetric tensor $S$
of order $p$, we define
\beqn{Tnorm}
\|S\| \eqdef \max_{\|v\|=1}  | S [v]^p |
= \max_{\|v_1\|= \cdots= \|v_p\|=1} | S[v_1, \ldots, v_p] |
\eeqn
the induced Euclidean norm (see \cite[Theorem~2.1]{ZhanLingQi12} for a proof
of the second equality). We denote by $\nabla_x^\ell f(x)$ the $\ell$-th
order derivative of $f$ evaluated at $x$, noting that such a tensor is
always symmetric for any $\ell\geq 2$. The notation $\nabla_x^\ell f(x)[s]^\ell$
denotes this $\ell$-th derivative tensor applied to $\ell$ copies of the vector $s$.
All inexact quantities are indicated by an overbar. 
For a symmetric matrices $M$, $\lambda_{\min}(M)$ is the smallest eigenvalue of
$M$. We will also use the function
\beqn{chi-def}
\chi_j(t) \eqdef \sum_{\ell = 1}^j \frac{t^\ell}{\ell !} \ms (t\geq 0),
\eeqn
where $j \geq 1$. We use the notation $\expect[X]$ to indicate the expected
value of a random variable $X$. In addition, given a random event $A$, $\prob(A)$
denotes the probability of $A$, while $\indic{A}$ refers to the indicator of the
random event $A$ occurring. The notation $\no{A}$ indicates that event $A$
does not occur.

\numsection{A regularization algorithm with inexact evaluations}\label{algo-s}

\subsection{The problem class}

We first make our framework more formal by detailing our assumptions on
problem \req{problem}.

\begin{description}
\item[AS.1] The function $f$ is $p$ times continuously
  differentiable in $\Re^n$. Moreover, its $j$-th
  order derivative tensor is  Lipschitz continuous for $j\in\ii{p}$ in the
  sense that there exist constants $L_{f,\ell}\geq 0$  such that, for all
  $\ell\in\ii{p}$ and all $x,y \in \Re^n$,
  \beqn{f-holder}
  \|\nabla_x^\ell f(x) - \nabla_x^\ell f(y)\| \leq L_{f,j} \|x-y\|.
  \eeqn
\item[AS.2] $f$ is bounded below in $\Re^n$, that is there exists a constant
  $\flow$ such that $f(x)\geq \flow$ for all $x \in \Re^n$.
\end{description}

\noindent
Because of AS.1, the $\ell$-th derivative of $f$ exists for $\ell\in\ii{p}$
and is a symmetric tensor of dimension $\ell$, which we denote by
\[
\nabla_x^\ell f(x) \eqdef \left(\frac{\partial^\ell f}{\partial x_{i_1}\ldots
  \partial x_{i_\ell}}\right)_{i_j\in\ii{n},j\in\ii{\ell}}(x).
\]
Moreover, the $p$-th degree Taylor series of $f$ at a point $x$ and
evaluated for a step $s$ is well-defined and can be written as
\beqn{Tf-def}
T_{f,p}(x,s)
\eqdef f(x)+\sum_{\ell= 1}^p\frac{1}{\ell!}\nabla_x^\ell f(x)[s]^\ell,
\eeqn
where $\nabla_x^\ell f(x)[s]^\ell$ denotes the scalar obtained by applying the
$\ell$-dimensional tensor $\nabla_x^\ell f(x)$ to $\ell$ copies of the vector
$s$. Because we will reuse this notation later, note that the first subscript
in $T_{f,p}(x,s)$ is the function whose Taylor expansion is
being considered, while the second is the degree of the expansion. The argument $x$ is
the point at which derivatives of $f$ are computed and $s$ is a step from $x$
so that $T_{f,p}(x,s)$ approximates the value of $f$ at the point $x+s$.
We will also make frequent use of the \emph{Taylor decrement} defined as
\beqn{DT-def}
\DT_{f,p}(x,s)
= T_{f,p}(x,0) - T_{f,p}(x,s)
=  - \sum_{\ell=1}^p \frac{1}{\ell!}\nabla_x^\ell f(x)[s]^\ell
\eeqn
We will also rely on the following well-known but important consequence of
AS.1.

\llem{Lipschitz-bounds-l}{
  Suppose that AS.1 holds. Then, for all $x,s \in \Re^n$,
  \beqn{f-Lip-bound}
  f(x+s) \leq T_{f,p}(x,s) + \frac{L_{f,p}}{(p+1)!}\|s\|^{p+1}
  \eeqn
  and
  \beqn{ders-Lip-bound}
  \|\nabla_x^\ell f(x+s) - \nabla_s^\ell T_{f,p}(x,s)\|
  \leq \frac{L_{f,p}}{(p-\ell+1)!}\|s\|^{p-\ell + 1}
  \tim{for all} \ell \in \ii{p}.
  \eeqn
}

\proof{See \cite[Lemma~2.1]{CartGoulToin20b}.} 

\subsection{Optimality measures}

We now turn to the important question of defining what we mean by
(approximate) critical points of arbitrary order but first address the
motivation for considering this issue.  In the standard exact case, it has
long been known that using Newton's method (i.e.\ a model of degree two)
practically outperforms the steepest descent method (which only uses a model of
degree one), even for computing first-order approximate critical points.  More
recently, it was shown in \cite{BirgGardMartSantToin17} that using a model of
degree $p>2$ (if possible) results in further improvements in evaluation
complexity. However, if an algorithm uses a model of degree $p>2$, why should it
be constrained to seek only for first- or second-order approximate critical
points? As it turns out, this question raises a number of issues, the first
being to define what is meant by an approximate critical point of general
order $q\leq p$. In the rest of this paper, we use the concept of
approximate minimizers discussed in \cite{CartGoulToin20a}.  Specifically,
given ``accuracy requests" $\epsilon= (\epsilon_1,\ldots, \epsilon_q)$ and
``optimality radii'' $\delta=(\delta_1,\ldots, \delta_q)$ with
\[
\epsilon_j \in (0,1] \tim{and} \delta_j \in (0,1] \tim{for} j \in \ii{q},
\]
we say that $x$ is a $q$-th order $(\epsilon,\delta)$-approximate
minimizer (or $(\epsilon,\delta)$-approximate
critical point) for problem \req{problem} if 
\beqn{strong}
\phi_{f,j}^{\delta_j} (x) \leq \epsilon_j \frac{\delta_j^j}{j!}
\tim{ for } j \in \ii{q},
\eeqn
where
\beqn{phi-def}
\phi_{f,j}^{\delta_j}(x) \eqdef
f(x) - \min_{\|d\| \leq \delta_j} T_{f,j}(x,d)
= \max_{\|d\| \leq \delta_j} \DT_{f,j}(x,d),
\eeqn
where, as is standard, the min and max are considered global.
Note that $\phi_{f,j}^\delta(x)$ is nothing but the largest decrease
obtainable on the $j$-th degree Taylor expansion of $f$ in a neighbourhood of
size $\delta_j$. As such, it is always well-defined for functions satisfying
AS.1 and is always non-negative. Also note that, because of the Cauchy-Schwarz
inequality, 
\beqn{strong1}
\phi_{f,1}^{\delta_1}(x)
= \max_{\|d\| \leq \delta_1} (-\nabla_x^1f(x)[d] )
= \|\nabla_x^1f(x)\| \,\delta_1,
\eeqn
and we immediately see that \req{strong}, when specialized to first-order, is
identical to the classical condition asking that $\|\nabla_x^1f(x)\|\leq
\epsilon_1$. Similarly, it is easy to verify that, when $\nabla_x^1f(x)=0$,
the second-order version of \req{strong} gives that
\[
\phi_{f,2}^{\delta_2}(x)
= \max_{\|d\| \leq \delta_2} (-\half \nabla_x^2f(x)[d]^2)
= \half \max\big[0, -\lambda_1[\nabla_x^2f(x)]\big]\,\delta_2^2
\]
where $\lambda_1[\nabla_x^2f(x)]$ is the leftmost eigenvalue of the Hessian
$\nabla_x^2f(x)$, so we obtain that \req{strong} is the same as the classical
condition that the absolute value of this eigenvalue is less than $\epsilon_2$
in this case. For example, the origin is both an $(\epsilon_1,1)$-approximate
first-order and an $(\epsilon_2,1)$-approximate
second-order minimizer of function $\sfrac{1}{6} x^3$ for any $\epsilon_1,\epsilon_1\in
(0,1]$, but is not an $(\epsilon_3,\delta_3)$-approximate
first-order one for any $\epsilon_3, \delta_3 \in (0,1]$.
We refer the reader to \cite{CartGoulToin20a} for a more
extensive discussion.

The condition \req{strong} has clear advantages over the  more usual definitions for first- and
second-order approximate critical points: it is well defined 
for all orders and it is a continuous\footnote{Difficulties
with the standard definition already start with order three because the
nullspace of $\nabla_x^2f(x)$ is not a continuous function of $x$.} function
of $x$.  Moreover, its evaluation is straightforward for $j=1$ (see
\req{strong1}) and easy for $q=2$ (it then reduces to the standard
trust-region subproblem whose cost is comparable to that of computing the
Hessian's leftmost eignevalue, see \cite[Chapter~7]{ConnGoulToin00}).
However, its evaluation may actually be extremely costly for $j>2$. From a
formal point of view, this does not affect the evaluation complexity of an
algorithm using it because it does not involve any new evaluation of $f$ and
its derivatives. We also note that we could consider an approximate version of
\req{strong}, where we would require that, for each $j\in \ii{q}$,  there
exists a $d_j$ such that $\|d_j\| \leq \delta_j$ and
\beqn{strong-approx}
\nu \phi_{f,j}^{\delta_j}(x)
\leq \DT_{f,j}(x,d_k)
\leq \epsilon_j \frac{\delta_j^j}{j!},
\eeqn
where $\nu$ is a constant in $(0,1]$. Note that \req{strong-approx} does not assume the
knowledge of the global minimizer or $\phi_{f,j}^{\delta}(x)$, but merely
that we can ensure the second part of \req{strong-approx} (see
\cite{deKlLaur19,deKlLaur20,SlotLaur20} for research in this direction). Note
also that, by definition, 
\[
\DT_{f,j}(x,d_j) \leq \nu \epsilon_j \frac{\delta_j^j}{j!}
\tim{ implies }
\phi_{f,j}^{\delta}(x) \leq \epsilon_j \frac{\delta_j^j}{j!}
\]
and this approximate and potentially less costly variant of \req{strong}
could thus replace it at the price of multiplying every $\epsilon_j$ by
the constant $\nu$. We will however ignore this possibility in our analysis,
keeping \req{strong} for simplicity of exposition.

\subsection{The regularization algorithm}

We are now in a position to describe our adaptive regularization
algorithm \algn whose purpose is to compute a $q$-th order
$(\epsilon,\delta)$-approximate minimizer of $f$ in problem~\req{problem}.
The vector of accuracies $\epsilon$ is given, together with a model degree $p\geq q$, corresponding to the maximum
order of available derivatives. If values of the objective function $f$ and its
derivatives of orders ranging from one to $p$ were known exactly, a typical
adaptive regularization method could be outlined as follows. At iteration
$k$, a local model of the objective function's variation would first be
defined by regularizing the Taylor series of degree $p$ at the current iterate
$x_k$, namely
\beqn{model-exact}
m_k(s) = -\DT_{f,p}(x_k,s) + \frac{\sigma_k}{(p+1)!}\|s\|^{p+1},
\eeqn
where $\sigma_k$ is a regularization parameter to be specified later.
A step $s_k$ would next be computed by approximately minimizing $m_k(s)$
in the sense that $m_k(s_k)\leq m_k(0) = 0$ and
\beqn{step-term-exact}
\phi_{m_k,j}^{\delta_{k,j}}(s_k) \leq \theta \epsilon_j \frac{\delta_{k,j}^j}{j!},
\eeqn
for some $\theta \in (0,\half)$ and $\delta_k \in (0,1]^q$. In this condition,
\[
\phi_{m_k,j}^{\delta_{k,j}}(s_k) = \max_{\|d\| \le \delta_{k,j}} \Delta T_{m_{k,j}}(s_k,d)
\]
is the $j$-th order optimality measure
\req{phi-def} for the model \req{model-exact} computed at $s_k$,
in which, for $j\in\ii{q}$, 
\beqn{Tm-def}
T_{m_k,j}(s_k,d)= m_k(s_k)
+ \sum_{\ell=1}^j\frac{1}{\ell!}\nabla_s^\ell T_{f,p}(x_k,s_k)[d]^\ell
+ \frac{\sigma_k}{(p+1)!} \sum_{\ell=1}^j\frac{1}{\ell!} \nabla_s^\ell \|s_k\|^{p+1}[d]^\ell
\eeqn
{for $d\in\Re^n$} and thus
\[
\DT_{m_k,j}(s_k,d)= - \sum_{\ell=1}^j\frac{1}{\ell!}\nabla_s^\ell T_{f,p}(x_k,s_k)[d]^\ell
- \frac{\sigma_k}{(p+1)!} \sum_{\ell=1}^j\frac{1}{\ell!} \nabla_s^\ell
\|s_k\|^{p+1}[d]^\ell
\]
(note the reuse of the notations introduced in \req{Tf-def} and \req{DT-def}, but for the
function $m_k(s)$ instead of $f(x)$).
The values of
$f(x_k+s_k)$ and $\{\nabla_x^\ell f(x_k+s_k)\}_{\ell=q+1}^p$ would then be
computed and the trial point $x_k+s_k$ would then be accepted as the next
iterate, provided the ratio
\[
\rho_k= \frac{f(x_k)-f(x_k+s_k)}{\DT_{f,p}(x_k,s_k)},
\]
is sufficiently positive. The regularization parameter $\sigma_k$ would
then be adapted/updated before a new iteration is started, providing the
``adaptive regularization'' suggested by the name of the method. (See
\cite{CartGoulToin20b} for the complete description of such an algorithm using
exact function and derivatives values.) The \algn algorithm follows the same lines, except that the values
of $f(x_k)$, $f(x_k+s_k)$ and $\DT_{f,p}(x_k,s_k)$ are not known 
exactly, the inexactness in the latter resulting from the inexactness of the
derivatives $\{\nabla_x^\ell f(x_k)\}_{\ell=1}^p$ .  
Instead, inexact values
$\barf(x_k)$, $\barf(x_k+s_k)$ and $\barDT_{f,p}(x_k,s_k)$ are now computed
and used to (re)-define the model
\beqn{model}
m_k(s) = -\barDT_{f,p}(x_k,s) + \frac{\sigma_k}{(p+1)!}\|s\|^{p+1}.
\eeqn
In particular, setting
 \beqn{omega-def}
0< \omega < \min \left[\frac{1-\eta}{3},\frac{\eta}{2}\right],
 \vspace*{-2mm}
 \eeqn
the approximations $\barf(x_k)$ and $\barf(x_k+s_k)$ are required to satisfy
the accuracy conditions
\begin{eqnarray}
\left|\overline{f}(x_k)-f(x_k)\right|&\leq& \omega\overline{\DT}_{f,p}(x_k,s_k),\label{f-est1}\\
\left|\overline{f}(x_k+s_k)-f(x_k+s_k)\right|&\leq& \omega\overline{\DT}_{f,p}(x_k,s_k).\label{f-est2}
\end{eqnarray}
In what follows, we will consistently denote inexact values by an overbar.

The model \req{model} is then approximately minimized by the feasible step
$s_k$ in the sense that the trial point $x_k+s_k$ satisfies
\beqn{descent}
m_k(s_k) \leq  m_k(0)=0
\eeqn
and
\beqn{step-term}
\barphi_{m_k,j}^{\delta_{k,j}}(s_k)
= \max_{\|d\|\leq \delta_{k,j}}\barDT_{m_k,j}(s_k,d)
\leq \theta \epsilon_j \frac{\delta_{k,j}^j}{j!},
\eeqn
for $j\in \ii{q}$ and some $\theta\in (0,\half)$ and $\delta_k \in (0,1]^q$.
The values $\barf(x_k)$, $\barf(x_k+s_k)$ and $\barDT_{f,p}(x_k,s_k)$ are also used to
compute the ratio $\rho_k$, the value of which decides of the acceptance of
the trial point.
The \algn algorithm is detailed as Algorithm~\ref{algnbr} \vpageref{algnbr}.

\algo{algnbr}{The \algn Algorithm}{
\vspace*{-4mm}
\begin{description}
\item[Step 0: Initialization.]
  An initial point $x_0\in\Re^n$, an initial regularization parameter
  $\sigma_0>0$ and a sought optimality order $q\in\ii{p}$ are given, as well as a vector of accuracies $\epsilon \in
  (0,1]^q$. The constants $\theta\in(0,\half)$, $\eta \in (0,1)$,  $\gamma>1$,
  $\alpha \in (0,1)$, {{$0< \omega < \min \left[\frac{1-\eta}{3},\frac{\eta}{2}\right]$}}
  and $\sigma_{\min}\in (0,\sigma_0)$ are also given. Set $k=0$.
\item[Step 1: Model construction. ] Compute approximate derivatives \\
  $\{\overline{ \nabla_x^\ell f}(x_k)\}_{\ell\in\{1,...,p\}}$ and form the model
  $m_k(s)$ defined in \eqref{model}.
\item[Step 2: Step calculation. ]
  Compute a step $s_k$ satisfying \req{descent} and
  \req{step-term} for $j\in\ii{q}$ and some $\delta_k \in (0,1]^q$. If
    $\barDT_{f,p}(x_k,s_k) = 0$, go to Step~4.
\item[Step 3: Function estimates computation. ] 
  Compute the approximations $\barf(x_k)$ and $\barf(x_k+s_k)$ of
  $f(x_k)$ and $f(x_k+s_k)$, respectively, such that
  \eqref{f-est1}--\eqref{f-est2} are satisfied. 
\item[Step 4: Acceptance test. ] Set
  \beqn{rho-def}
  \rho_k = \left\{\begin{array}{ll}
  \bigfrac{\barf(x_k) - \barf(x_k+s_k)}{\barDT_{f,p}(x_k,s_k)} &
  \tim{ if } \barDT_{f,p}(x_k,s_k) > 0, \\
  -\infty & \tim{ otherwise.}
  \end{array}\right.
  \eeqn
  If $\rho_k \geq \eta$ (\textit{successful iteration}), then define
  $x_{k+1} = x_k + s_k$; otherwise (\textit{unsuccessful iteration}) define $x_{k+1} = x_k$.
\item[Step 5: Regularization parameter update. ]
  Set
 \vspace*{-2mm}
 \beqn{sigupdate}
 \sigma_{k+1} = \left\{ \begin{array}{ll}
   \max\left[\sigma_{\min},\bigfrac{1}{\gamma} \sigma_k\right],     &\tim{if} \rho_k \geq \eta,\\
   \gamma \sigma_k, & \tim{if} \rho_k < \eta.
 \end{array} \right.
 \eeqn
 
{{ Increment $k$ by one and go to Step~1.}}
\end{description}
}

\noindent We first verify that the algorithm is well-defined.

\llem{step-exists}{A step $s_k$ satisfying \req{descent} and
  \req{step-term} for $j\in\ii{q}$ and some $\delta_k \in (0,1]^q$ always exists.   
}
  \proof{The proof is a direct extension of that of \cite[Lemma~4.4]{CartGoulToin20a} using
    inexact models. It is given in appendix for completeness. 
  } 
  
\noindent
Some comments on this algorithm are useful at this stage.
\begin{enumerate}
\item
It is important to observe that the algorithm is fully implementable with
existing computational technology in the very frequent cases where $q=1$ or
$q=2$. Indeed the value of $\barphi_{m_k,1}^{\delta_{k,1}}$ can easily be
obtained analytically.  When $q=2$, the same comment obviously applies for
$\barphi_{m_k,1}^{\delta_{k,1}}$, while the value $\barphi_{m_k,2}^{\delta_{k,2}}$
can be computed by a standard trust-region solver (whose cost is comparable
to that of the more usual calculation of the most negative eigenvalue), again
making the algorithm practical.
We refer the interested reader to \cite{BellGuriMoriToin21} for the
presentation of numerical results in the framework of finite sum optimization
for automatic learning.

In other cases, the computation $\barphi_{m_k,j}^{\delta_{k,j}}$ may be
extremely expensive, making our approach mostly 
theoretical at this stage.  However, we recall that, since evaluations of the
objective function and its derivatives do not occur in this computation
(once the approximate derivatives are known), its cost has
no impact on the evaluation complexity of the \algn algorithm.

\item
We assume in what follows that, once the inexact model $m_k(s)$ is
determined, then the computation of the pair $(s_k,\delta_k)$ (and thus of the 
trial point $x_k+s_k$) is deterministic.  Moreover, we assume that the
mechanism which ensures \req{f-est1}-\req{f-est2} in Step~3 of the algorithm is also
deterministic, so that $\rho_k$ and the fact that iteration $k$ is successful
are deterministic outcomes of the realization of the inexact model.  

\item
Observe that, because we have chosen $m_k$ to be a model of the local
variation in $f$ rather than a model of $f$ itself, $\barf(x_k)$ is
not needed (and not computed) in Steps~1 and 2 of the algorithm.  This
distinguishes the \algn algorithm from the approaches of
\cite{BlanCartMeniSche19,ChenMeniSche18}.
\end{enumerate}
In what follows, all random quantities are denoted by capital letters, while
the use of small letters is reserved for their realization.  In particular,
let us denote a random model at iteration $k$ as $M_k$, while we use the
notation $m_k$ for its realizations. Given $x_k$, the source of randomness in
$m_k$ comes from the random approximation of the derivatives.  Similarly, the
iterates $X_k$, as well as the regularization parameters $\Sigma_k$ and the
steps $S_k$ are random variables (except for initial values $x_0$ and
$\sigma_0$ for the former two) and $x_k$, $\sigma_k$ and $s_k$ denote their
realizations. Moreover, $\delta_k$ denotes a realization of the random vector
$\Delta_k$ arising in (\ref{step-term}).  Hence, the \algn Algorithm generates
a random process
\begin{eqnarray}
  {{\{X_k, S_k, M_k, \Sigma_k,\Delta_k
  \}.\label{sprocess} }}
\end{eqnarray}
where $X_0 = x_0$ and $\Sigma_0= \sigma_0$ are deterministic.

\subsection{The probabilistic setting}

We now make our probabilistic assumptions on the
\algn algorithm explicit. For $k\geq 0$, our assumption on the past is
formalized by considering $\sigax$ the $\hat{\sigma}$-algebra
induced by the random variables $M_0$, $M_1$,..., $M_{k-1}$, with
$\calA_{-1}^{M}=\hat{\sigma}(x_0)$.
In order to formalize our probabilistic assumptions we need a few more definitions.
We define, at iteration $k$ of an
arbitrary realization,
\beqn{dkj-def}
d_{k,j} = \argmax_{\|d\|\leq\delta_{k,j}}\DT_{m_k,j}(s_k,d)
\eeqn
the argument of the maximum in the definition of
$\phi_{m_k,j}^{\delta_{k,j}}(x_k)$, and
\beqn{bardkj-def}
\bard_{k,j} = \argmax_{\|d\|\leq\delta_{k,j}}\barDT_{m_k,j}(s_k,d)
\eeqn
that in the definition of $\barphi_{m_k,j}^{\delta_{k,j}}(s_k)$.
We also define, at the end of Step 2 of iteration $k$, the events
\begin{equation}
   \label{defAMk}
   \AM_k= \left\{\begin{array}{ll}
   \AM_k^{(1)}\cap \bigcap_{j=1}^q\left(\AM_{k,j}^{(2)}\cap\AM_{k,j}^{(3)}\right)
   & \tim{if $q\in \{1,2\}$}\\
   \AM_k^{(1)}\cap \AM_k^{(4)}\cap\bigcap_{j=1}^q\left(\AM_{k,j}^{(2)}\cap\AM_{k,j}^{(3)}\right)
   &\tim{ otherwise,}
   \end{array}\right.
\end{equation}
with
\[
\begin{array}{lcl}
\AM_k^{(1)}&=&\left\{|\barDT_{f,p}(X_k,S_k)-\DT_{f,p}(X_k,S_k)|
\leq \omega \barDT_{f,p}(X_k,S_k)\right\},\\
\AM_{k,j}^{(2)}&=&\big\{|\barDT_{m_k,j}(S_k,D_{k,j})-\DT_{m_k,j}(S_k,D_{k,j})|
  \leq \omega \barDT_{m_k,j}(S_k,D_{k,j}),\\
\AM_{k,j}^{(3)}&=&\big\{|\barDT_{m_k,j}(S_k,\barD_{k,j})-\DT_{m_k,j}(S_k,\barD_{k,j})|
   \leq \omega \barDT_{m_k,j}(S_k,\barD_{k,j}),\\
\AM_{k}^{(4)}&=&\big\{  \max_{\ell \in \iibe{2}{p}} \|\overline{\nabla_x^\ell f}(X_k) \| \leq \Theta  \},
\end{array}
\]
for some $ \Theta>0 $.  Note that  $ \Theta $ is independent of $k$ and
does not need to be known explicitly. Moreover, $\AM_k^{(4)}$ is not involved in the definition of $\AM_k$
if $q\in \{1,2\}$. In what follows, we will say that
iteration $k$ is \textit{accurate}, if $\indicM=1$, and
iteration $k$ is \textit{inaccurate}, if $\indicM=0$.

The conditions defining $\AM_k$ may seem abstract at first sight, but we now motivate them by
looking at what kind of accuracy on each derivative
$\overline{\nabla_x^\ell f}(x_k)$  ensures that they hold.

\llem{suff-acc-conds}{For each $k\geq 0$, we have the following.
  \begin{enumerate}
  \item Let
  \beqn{tau-def}
  \tau_k \eqdef \max\left[ \|S_k\|, \max_{j\in\ii{q}}[\|D_{k,j}\|, \|\barD_{k,j}\|]\right]
  \eeqn
  and
  \beqn{DTstar-def}
  \barDT_{k,\min}
  \eqdef \min\left[\barDT_{f,p}(X_k,S_k),\min_{j\in\ii{q}}\Big[ \barDT_{m_k,j}(S_k,D_{k,j}),
                   \barDT_{m_k,j}(S_k,\barD_{k,j})\Big]\right].
  \eeqn
  Then $\AM_k^{(1)}$, $\{\AM_{k,j}^{(2)}\}_{j=1}^q$ and
  $\{\AM_{k,j}^{(3)}\}_{j=1}^q$ occur if
  \beqn{the-cond}
  \|\overline{\nabla_x^\ell f}(X_k) -\nabla_x^\ell f(X_k)\|
  \leq \omega \frac{\barDT_{k,\min}}{6\tau_k^\ell}
  \tim{ for }\ell\in\ii{p}.
  \eeqn
\item Suppose that AS.1 holds.
  Then $\AM_k^{(4)}$ occurs if
    \beqn{conds-derm}
    \|\overline{\nabla_x^\ell f}(X_k) -\nabla_x^\ell f(X_k)\|
    \leq \Theta_0
    \tim{ for }\ell\in\iibe{2}{p}
    \eeqn
    and some constant $\Theta_0\geq 0$ independent of $k$ and $\ell$.
  \end{enumerate}
}

\proof{Consider the first assertion.  That $\AM_k^{(1)}$ occurs follows from the inequalities
\[
\begin{array}{lcl}
|\barDT_{f,p}(X_k,S_k)-\DT_{f,p}(X_k,S_k)|
& \leq & \bigsum_{\ell=1}^p \bigfrac{\|S_k\|^\ell}{\ell!}\,
         \|\overline{\nabla_x^\ell f}(X_k) - \nabla_x^\ell f(X_k)\|\\*[2ex]
& \leq & \bigsum_{\ell=1}^p \bigfrac{\tau_k^\ell}{\ell!}\,
         \|\overline{\nabla_x^\ell f}(X_k) - \nabla_x^\ell f(X_k)\|\\*[2ex]
& \leq & \bigsum_{\ell=1}^p \bigfrac{\omega}{6\ell!}\,
         \barDT_{k,\min}\\*[2ex]
& \leq & \bigsum_{\ell=1}^p \bigfrac{\omega}{6\ell!}\,
         \barDT_{f,p}(X_k,S_k)\\*[2ex]
& \leq & \bigfrac{1}{6}\chi_p(1)\,\omega\,\barDT_{f,p}(X_k,S_k)\\*[2ex]
&   <  & \omega \,\barDT_{f,p}(X_k,S_k).
\end{array}
\]
where we have used \req{tau-def}, \req{the-cond}, \req{DTstar-def} and the
fact that $\chi_p(1)\leq 2$. The verification that $\{\AM_{k,j}^{(2)}\}_{j=1}^q$ and
$\{\AM_{k,j}^{(3)}\}_{j=1}^q$ also occur uses a very similar
argument, with one additional ingredient: employing the triangle inequality,
\req{model}, we have that, for all $\ell\in\ii{p}$,
\[
\left\|\overline{\nabla_d^\ell T}{\,}_{m_k,j}(S_k,0)-\nabla_d^\ell T_{m_k,j}(S_k,0)\right\|
\leq \bigsum_{t=\ell}^p
  \left\|\overline{\nabla_x^t f}(X_k)-\nabla_x^t f(X_k)\right\|
         \bigfrac{\|S_k\|^{t-\ell}}{(t-\ell)!}.
\]
Considering now $D =D_{k,j}$ or $D=\barD_{k,j}$ and using the above inequality,
\req{tau-def}, \req{the-cond}, \req{DTstar-def} and the
facts that $\chi_j(1)\leq 2$ and $\chi_{p-\ell}(1)\leq 2$, we have that
\[
\begin{array}{l}
|\barDT_{m_k,j}(S_k,D)-\DT_{m_k,j}(S_k,D)|\\*[2ex]
\hspace*{30mm} \leq \bigsum_{\ell=1}^j \bigfrac{\|D\|^\ell}{\ell!}
        \|\overline{\nabla_d^\ell T_{m_k,j}}(S_k,0) - \nabla_d^\ell T_{m_k,j}(S_k,0)\|\\*[2ex]
\hspace*{30mm} \leq \bigsum_{\ell=1}^j \bigfrac{\|D\|^\ell}{\ell!}
        \bigsum_{t=\ell}^p
  \left\|\overline{\nabla_x^t f}(X_k)-\nabla_x^t f(X_k)\right\|
         \bigfrac{\|S_k\|^{t-\ell}}{(t-\ell)!}\\[2ex]
\hspace*{30mm} \leq \bigsum_{\ell=1}^j \bigfrac{1}{\ell!}
        \bigsum_{t=\ell}^p
  \left\|\overline{\nabla_x^t f}(X_k)-\nabla_x^t f(X_k)\right\|
         \bigfrac{\tau_k^t}{(t-\ell)!}\\[2ex]
\hspace*{30mm} \leq \bigsum_{\ell=1}^j \bigfrac{1}{\ell!}
        \bigsum_{t=\ell}^p\bigfrac{1}{(t-\ell)!}\,
        \omega\frac{\barDT_{k,\min}}{6}\\*[2ex]
\hspace*{30mm} \leq \bigfrac{1}{6}\omega\,\barDT_{k,\min}
        \bigsum_{\ell=1}^j \bigfrac{1}{\ell!} (1+\chi_{p-\ell}(1))\\*[2ex]
\hspace*{30mm} \leq \omega \, \barDT_{m_k,j}(S_k,D),
\end{array}
\]
as desired. To prove the second assertion, observe that AS.1
implies that $\|\nabla_x^\ell f(X_k)\| \leq L_{f,\ell-1}$  for $j\in \iibe{2}{p}$,
and thus, using \req{conds-derm}, that,
for $\ell \in \iibe{2}{p}$,
\[
\begin{array}{lcl}
\|\overline{\nabla_x^\ell f}(X_k) \|
& \leq & \|\nabla_x^\ell f(X_k)\|+\|\overline{\nabla_x^\ell f}(X_k)-\nabla_x^\ell f(X_k)\|\\*[1.5ex]
& \leq & L_{f,\ell-1} +\Theta_0.
\end{array}
\]
This gives the desired conclusion with the choice $\Theta =
\max_{\ell\in\iibe{2}{p}}L_{f,\ell-1}+ \Theta_0$.
} 

\noindent
Of course, the conditions stated in Lemma~\ref{suff-acc-conds} are sufficient
but by no means necessary to ensure $\AM_k$.  In particular, they make no
attempt to exploit a possible favourable balance between the errors made on
derivatives at different degrees, nor do they take into account that
$\AM_k^{(1)}$, $\AM_{k,j}^{(2)}$ and $\AM_{k,j}^{(3)}$ only specify conditions
on model accuracy in a finite, dimension-independent subset of
directions. Despite these limitations, \req{the-cond} and \req{conds-derm}
allow the crucial conclusion that $\AM_k$ does occur if the derivatives
$\overline{\nabla_x^jf}(X_k)$ are sufficiently accurate compared to the model
decrease.  Moreover, since one would expect that, as an approximate minimizer
is approached, $\|S_k\|$, $\|D_{k,j}\|$ and $\|\barD_{k,j}\|$ (and thus
$\tau_k$) become small, they also show the accuracy requirement becomes looser
for derivatives of higher degree.

\noindent
We now formalize our assumption on the stochastic process generated by the
\algn algorithm.

\begin{description}
\item[AS.3]\mbox{}\\
For all $k\ge 0$, the event $ \AM_k$ satisfies the condition
  \beqn{prob-conds}
   p_{\calM,k} = \prob( \AM_k | \calA_{k-1}^{M} )=\expect[\indicM|\calA_{k-1}^{M}]\geq \pM
   \eeqn
   for some $\pM \in (\half,1]$ independent of $k$.
 \end{description}

\noindent
AS.3 is realistic for instance in cases where derivatives are approximated by
randomized finite differences or by subsampling in the context of finite sum
minimization.

We observe that, in contrast with \cite{BlanCartMeniSche19,ChenMeniSche18}, the definition of $\AM_k$  does not require the
model to be ``linearly/quadratically'' accurate everywhere in a ball around
$x_k$ of radius at least $\|s_k\|$, but merely that their variation is
accurate enough along $s_k$ (as specified in $\AM_k^{(1)}$) and along
$d_{k,j}$ and $\bard_{k,j}$ (as specified in $\AM_{k,j}^{(2)}$ and
$\AM_{k,j}^{(3)}$ )\footnote{A slightly stronger assumption would be to 
require a sufficient relative accuracy along $s_k$ and in
a (typically small) neighbourhood of $s_k$.} for all $j\in \ii{q}$.
The need to consider $\AM_{k,j}^{(2)}$ and $\AM_{k,j}^{(3)}$ for $j\in \ii{q}$
in the definition of $\AM_k$ results from our insistence that $q$-th order
approximate optimality must include $j$-th order approximate optimality for
all such $j$. AS.3 also parallels assumptions in
\cite{BlanCartMeniSche19,CartSche17,ChenMeniSche18,PaquSche20} where accuracy 
in derivatives' values is measured using the guaranteed model decrease or proxies
given by the $(p+1)$-st power of the trust-region radius or
the steplength. Finally, the conditions imposed by $\AM_{k,j}^{(2)}$ and
$\AM_{k,j}^{(3)}$ are only used whenever considering the value of
$\barphi_{m_k,j}^{\delta_{k,j}}(s_k)$, that is in Lemma~\ref{phim-barphim-l},
itself only called upon in Lemma~\ref{step-inequality} in the case where
$\|S_k\| \leq 1$. As a consequence, they are irrelevant when long steps are taken
($\|S_k\|>1$).   

\numsection{Worst-case evaluation complexity}\label{complexity-s}

Having set the stage and stated our assumptions, we may now consider the
worst-case evaluation complexity of the \algn algorithm. Our aim is to derive
a bound on the expected number of iterations $\expect(N_\epsilon)$ which is
needed, in the worst-case, to reach an $(\epsilon,\delta)$-approximate
$q$-th-order-necessary minimizer. Specifically, $N_\epsilon$ is the number
of iterations required until \req{strong} holds for the first time, i.e.,
\beqn{Neps-def}
N_\epsilon
=\inf \left\{ k\geq 0~\mid~ \phi_{f,j}^{\Delta_{k-1,j}}(X_k) \leq \epsilon_j
\frac{\Delta_{k-1,j}^j}{j!} \tim{for} j \in \ii{q}\right\}.
\eeqn
Note that
$\phi_{f,j}^{\Delta_{k-1,j}}(X_k)$, the $j$-th order optimality measure at
iteration $k$, uses the optimality radii $\Delta_{k-1,j}$ resulting from
the step computation at iteration $k-1$, as is the case in
\cite{BellGuriMoriToin19,CartGoulToin20b}.
Now recall that the trial point $X_{k-1}+S_{k-1}$ and the vector of radii $\Delta_{k-1}$
are deterministic once the inexact model at iteration $k-1$ is known.  Thus
these variables are measurable for $\sigax$ and because of our
deterministic assumptions on the accuracy of $f$, the event
$\{X_k=X_{k-1}+S_{k-1}\}$ (which occur when iteration $k-1$ is successful) is
also measurable for $\sigax$. As a consequence and since 
$\phi_{f,j}^{\Delta_{k-1,j}}(X_k)$ uses exact derivatives of $f$, the event
$\{N_\epsilon = k\}$ is measurable with respect to $\sigax$. The
definition \req{Neps-def} can thus be viewed as that of a
family of $\epsilon$-dependent stopping times for the stochastic process
generated by the \algn algorithm (see, e.g., \cite[section $2.3$]{CartSche17}).

\subsection{General properties of the \algn algorithm}

We first consider properties of ``accurate'' iterations, in the sense that
$\AM_k$ occurs, and start with the relation between
$\phi_{m_k,j}^{\delta_{k,j}}(s_k)$ and its approximation. The next lemma is
inspired by Lemma~3.2 in \cite{BellGuriMoriToin19}, but significantly differs
in that it now requires considering both directions $d_{k,j}$ and $\bard_{k,j}$.

\llem{phim-barphim-l}{
Consider any realization of the algorithm and assume that $\AM_k$ occurs. Then,
for $j\in \ii{q}$,
\beqn{phi-barphi}
\big(1-\omega\big)\barphi_{m_k,j}^{\delta_{k,j}}(s_k)
\leq \phi_{m_k,j}^{\delta_{k,j}}(s_k)
\leq \big(1+\omega\big) \barphi_{m_k,j}^{\delta_{k,j}}(s_k)
\eeqn
}

\proof{
Let $j\in \ii{q}$. Consider $d_{k,j}$ defined in \req{dkj-def}. From
\req{step-term}, we have that 
\[
\begin{array}{lcl}
\DT_{m_k,j}(s_k,d_{k,j})
& \leq & \barDT_{m_k,j}(s_k,d_{k,j})+|\DT_{m_k,j}(s_k,d_{k,j})-\barDT_{m_k,j}(s_k,d_{k,j})|\\*[2ex]
& \leq & \big(1+\omega\big)\barDT_{m_k,j}(s_k,d_{k,j})\\*[2ex]
& \leq & \big(1+\omega\big)
  \bigmax_{\|d\|\leq\delta_{k,j}}\barDT_{m_k,j}(s_k,d)\\*[2ex]
& = & \big(1+\omega\big)\barDT_{m_k,j}(s_k,\bard_{k,j})
\end{array}
\]
where we used the fact that $\AM_k$ occurs to derive the second inequality  and
considered $\bard_{k,j}$ defined in \req{bardkj-def}.
Therefore
\[
\phi_{m_k,j}^{\delta_{k,j}}(s_k)
= \DT_{m_k,j}(s_k,d_{k,j})
\leq \big(1+\omega\big)\barphi_{m_k,j}^{\delta_{k,j}}(s_k).
\]
This proves the rightmost inequality of \req{phi-barphi}.
Similarly, using our
assumption that $\AM_k$ occurs, we obtain that
\[
\begin{array}{lcl}
\DT_{m_k,j}(s_k,\bard_{k,j})
& \geq & \barDT_{m_k,j}(s_k,\bard_{k,j})
         - |\DT_{m_k,j}(s_k,\bard_{k,j})- \barDT_{m_k,j}(s_k,\bard_{k,j})|\\*[2ex]
& \geq & \big(1-\omega\big) \barDT_{m_k,j}(s_k,\bard_{k,j})\\*[2ex]
\end{array}
\]
and hence, from \req{phi-def} and \req{step-term}, that
\[
\big(1-\omega\big) \barphi_{m_k,j}^{\delta_{k,j}}(s_k)
\leq \max_{\|d\|\leq\delta_{k,j}} \DT_{m_k,j}(s_k,d)
= \phi_{m_k,j}^{\delta_{k,j}}(s_k),
\]
which concludes the proof of \req{phi-barphi}.
} 

The next step is to adapt an important property of $\Delta_{k,j}$ in the
exact case to our inexact framework.

\llem{delta-lower-l}{
  Suppose that AS.1 holds.  Then, for any $j\in\ii{q}$,
\begin{enumerate}
  \item if $j\in\{1,2\}$, $\Delta_{k,j}$ can always be
    chosen equal to one;
  \item in the other cases, and assuming that $\AM_k$ occurs,
    then, either $\|s_k\|> 1$ or $\Delta_{k,j} \leq 1$ can be chosen such that
  \beqn{delta-lower}
  \Delta_{k,j} \geq \kappa_\delta(\sigma_k)\epsilon_j,
  \eeqn
  where $\kappa_\delta(\sigma)\in (0,1)$ is independent of $\epsilon$ and
  decreasing when $\sigma$ grows.
  \end{enumerate}
 }

\proof{The proof broadly follows the developments of \cite[Lemmas~4.3 and
   4.4]{CartGoulToin20a}, except that it now uses the model involving
   approximate derivatives and that $L_f$, the upper bound of the
   derivatives of $f$ at $x_k$ derived from AS.1 is now replaced by
   $\Theta$, as guaranteed by $\AM_k^{(4)}$. The details (including the reason
   for the dichotomy between the two cases in the lemma's statement) are provided in
   appendix.}

\noindent
In what follows, we will assume that, whenever $q>2$,
the \algn algorithm computes a pair $(s_k,\delta_k)$ such that, for each
$j\in\ii{q}$, $\delta_{k,j}$ is always within a fraction of its maximal value,
thereby ensuring \req{delta-lower}. We now prove a crucial inequality relating
the step length to the accuracy requirements.

\llem{step-inequality}{
Consider any realization of the algorithm. Assume that $\AM_k$ occurs, that
iteration $k$ is successful and that, for some $j\in\ii{q}$, \req{strong}
fails for $(x_{k+1},\delta_{k,j})$.  Then either $\|s_k\| > 1$ or 
\beqn{central}
(1-2\theta)\epsilon_j\bigfrac{\delta_{k,j}^j}{j!}
\leq \frac{L_{f,p}+ \sigma_k}{(p-q+1)!} \bigsum_{\ell=1}^j\bigfrac{\delta_{k,j}^\ell}{\ell!}\|s_k\|^{p-\ell+1}
\eeqn
}
\proof{[See \cite[Lemma~5.3]{CartGoulToin20a} for the composite unconstrained
Lipschitz continuous case.] Suppose that $\|s_k\| \leq 1$. Since \req{strong}
fails at $(x_{k+1},\delta_{k,j})$, we must have that
\beqn{phi-fails}
\phi_{f,j}^{\delta_{k,j}}(x_{k+1}) > \epsilon_j \frac{\delta_{k,j}^j}{j!} > 0
\eeqn
for some $j\in\ii{q}$. Define $d$ to be the argument of the minimum in the definition of
$\phi_{f,j}^{\delta_{k,j}}(x_{k+1})$.  Hence, 
\beqn{d-bound}
0< \|d\| \leq \delta_{k,j}.
\eeqn
Using \req{phi-fails}, \req{phi-def} and the triangle
inequality, we thus obtain that 
\beqn{part0}
\phi_{f,j}^{\delta_{k,j}}(x_{k+1})
=  \DT_{f,j}(x_{k+1},d)
\leq \left|\DT_{f,j}(x_{k+1},d)-\DT_{m_k,j}(s_k,d)\right|
         + \DT_{m_k,j}(s_k,d).
\eeqn
Recalling now from \cite[Lemma~2.4]{CartGoulToin20b}) that
\[
\|\nabla_s^\ell\|s_k\|^{p+1}\|
= \frac{(p+1)!}{(p-\ell+1)!}\|s_k\|^{p-\ell+1},
\]
we
may now use the fact that $x_{k+1}=x_k+s_k$ since iteration $k$ is successful,
\req{ders-Lip-bound} in Lemma~\ref{Lipschitz-bounds-l}, \req{Tm-def},
\req{d-bound} and the triangle inequality to obtain that  
\beqn{part1}
\begin{array}{lcl}
\left|\DT_{f,j}(x_{k+1},d)-\DT_{m_k,j}(s_k,d)\right|
& \leq & \bigsum_{\ell=1}^j \frac{\delta_{k,j}^\ell}{\ell!}
          \|\nabla_x^\ell f(x_{k+1})-\nabla_s^\ell T_{f,p}(x_k,s_k)\|\\*[2ex]
&      &\hspace*{1cm} +\bigfrac{\sigma_k}{(p+1)!}
          \bigsum_{\ell=1}^j\bigfrac{\delta_{k,j}^\ell}{\ell!}
          \|\nabla_s^\ell\|s_k\|^{p+1}\|\\*[2ex]
& \leq & \bigfrac{L_{f,p}+ \sigma_k}{(p-q+1)!}
          \bigsum_{\ell=1}^j\bigfrac{\delta_{k,j}^\ell}{\ell!}\|s_k\|^{p-\ell+1}
\end{array}
\eeqn
Moreover, using \req{step-term}, \req{phi-barphi} and the fact that
{{$\omega< 1$}} (see \req{omega-def}), we deduce that
\beqn{part2}
\DT_{m_k,j}(s_k,d)
\leq  \phi_{m_k,j}^{\delta_{k,j}}(s_k)
\leq \big(1+\omega \big) \barphi_{m_k,j}^{\delta_{k,j}}(s_k)
\leq 2\theta \epsilon_j\bigfrac{\delta_{k,j}^j}{j!}.
\eeqn
Substituting \req{part1} and
\req{part2} into \req{part0} and using \req{d-bound} and \req{phi-fails}, we
obtain \req{central}.
} 

\llem{longs-l}{
  Suppose that AS.1 holds and consider any realization of the
  algorithm. Suppose also that $\AM_k$ occurs, that iteration $k$ is
  successful and that, for some $j\in\ii{q}$, \req{strong} fails
  for $(x_{k+1},\delta_{k,j})$. Then
  \beqn{longs}
  \|s_k\|^{p+1} \geq \psi(\sigma_k) \epsilon_j^\varpi
  \eeqn
  where
  \beqn{pi-def}
  \varpi = \left\{\begin{array}{ll}
  \bigfrac{p+1}{p-q+1} 
  & \tim{ if $q\in\{1,2\}$},\\*[2ex]
  \bigfrac{q(p+1)}{p}
  & \tim{otherwise.}
  \end{array}\right.
  \eeqn
  and
  \beqn{psi-def}
  \psi(\sigma) = \left\{\begin{array}{ll}
   \min\left[ 1,\left(\bigfrac{(1-2\theta)(p-q+1)!}
     {q!(L_{f,p}+\sigma)}\right)^\varpi\right]
   & \tim{if $q\in\{1,2\}$},\\*[3ex]
   \min\left[ 1,\left(\bigfrac{(1-2\theta)(p-q+1)!\,\kappa_\delta(\sigma)^{q-1}}
           {q!(L_{f,p}+\sigma)}\right)^\varpi\right]
   & \tim{otherwise.}
                   \end{array}\right.
  \eeqn
}

\proof{[See \cite[Lemma~5.4]{CartGoulToin20a}.] If
$\|s_k\| > 1$, the conclusion immediately follows. Suppose therefore that
$\|s_k\| \leq 1$ and consider $j$ such that \req{central} holds.
Recalling the definition of $\chi_j$ in \eqref{chi-def}, \req{central} can be
rewritten as  
\beqn{ineq-ab}
\alpha_k \,\epsilon_j \, \delta_{k,j}^j
\leq \|s_k\|^{p+1} \chi_j\left(\frac{\delta_{k,j}}{\|s_k\|}\right)
\eeqn
where  we have set
\[
\alpha_k = \frac{(1-2\theta)(p-q+1)!}{q!(L_{f,p}+\sigma_k)}.
\]
In particular, since  $\chi_j(t) \leq 2 t^j$
for $t \geq 1$, we have that, when $\|s_k\| \leq \delta_{k,j}$,
\beqn{ineq-good}
\alpha_k \,\epsilon_j
\leq 2 \|s_k\|^{p+1}\left(\frac{1}{\|s_k\|}\right)^j
= 2 \|s_k\|^{p-j+1}.
\eeqn
Suppose first that $q \in\{1,2\}$.  Then, from our assumptions and
Lemma~\ref{delta-lower-l}, $\delta_{k,j}=1$ and $\|s_k\| \leq 1 = \delta_{k,j}$.  Thus
\req{ineq-good} yields the first case of \req{pi-def}--\req{psi-def}. Suppose
now that $q>2$.  Then our assumptions imply that \req{delta-lower}
holds. If $\|s_k\|\leq \delta_{k,j}$, we may again deduce from
\req{ineq-good} that the first case of \req{pi-def}--\req{psi-def} holds,
which implies, because $\kappa_{\delta}(\sigma)<1$ and $1/(p-j+1) \leq
j/p$
, that
the second case also holds.  Consider therefore the case where $\|s_k\| >
\delta_{k,j}$.  Then \req{ineq-ab} and the fact that $\chi_j(t) < 2t$ for
$t\in [0,1]$ give that
\[
\alpha_k \,\epsilon_j \, \delta_{k,j}^j
\leq 2 \|s_k\|^{p+1} \left(\frac{\delta_{k,j}}{\|s_k\|}\right),
\]
which, with \req{delta-lower}, implies the second case of
\req{pi-def}--\req{psi-def} as requested.
}  

\noindent
Note that $\psi(\sigma)$ is decreasing as a function of $\sigma$ in both cases
of \req{psi-def}. We now investigate the decrease of the exact objective
function values at successful iterations.

\llem{accurate-f-decrease-l}{
  Suppose that AS.1 holds and consider any realization of the algorithm.  Then 
  \beqn{taylor-decrease}
  \barDT_{f,p}(x_k,s_k)
  \geq \frac{\sigma_k}{(p+1)!}\|s_k\|^{p+1}
  \geq \frac{\sigma_{\min}}{(p+1)!}\|s_k\|^{p+1}\ge 0,
  \eeqn
  where $\sigma_{\min}$ is defined in Step~0 of the \algn algorithm.
  Moreover, if iteration $k$ is successful, then
  \beqn{acc-f-decrease-1}
  f(x_k)-f(x_{k+1})
  \geq \frac{{{(\eta-2\omega)}}\sigma_{\min}}{(p+1)!}\|s_k\|^{p+1}{{> 0}}.
  \eeqn
}

\proof{
The inequality \req{taylor-decrease} immediately follows from
\req{model}, \req{descent}, \req{sigupdate}.
Now the fact that iteration $k$ is successful, together with \req{omega-def}
and \eqref{f-est1}--\eqref{f-est2}, imply that
\[
\begin{array}{lcl}
f(x_k)-f(x_{k+1})
& \geq & \barf(x_k) - \barf(x_{k+1}) - 2\omega\barDT_{f,p}(x_k,s_k) \\*[2ex]
& \geq & \eta \barDT_{f,p}(x_k,s_k){{-2\omega \barDT_{f,p}(x_k,s_k),}}
\end{array}
\]
yielding \req{acc-f-decrease-1} using \req{taylor-decrease} and {{\eqref{omega-def}}}.
} 

\noindent
We finally conclude our analysis of ``accurate'' iterations by proving a
standard result in the analysis of adaptive regularization methods. A similar
version of this result was presented in \cite[Lemma~4.2]{BellGuriMoriToin19}
for the case where both function values and models are sufficiently accurate.

\llem{acc-sigma-success-l}{
Suppose that AS.1 holds and let $\beta>1$ be given. Then, for any realization of the algorithm, if iteration
$k$ is such that $\AM_k$ occurs and 
\beqn{large-sigma}
\sigma_k \geq \sigma_s \eqdef \max\left[\beta\sigma_0,\frac{L_{f,p}}{{{1-\eta-3\omega}}}\right],
\eeqn
then iteration $k$ is successful.
}

\proof{
  Suppose that \req{large-sigma} holds.  Thus, using successively
  \req{rho-def}, the triangle inequality, the fact that $\AM_k$ occurs,
  \req{f-Lip-bound}, \req{taylor-decrease}, \req{omega-def}, \eqref{f-est1}--\eqref{f-est2} and
  \req{large-sigma}, we deduce that 
  \[
  \begin{array}{lcl}
  |\rho_k-1|
  & \leq & \bigfrac{1}{\barDT_{f,p}(x_k,s_k)}
           \Big[\Big(\barf(x_k) - f(x_k)\Big)
             +\Big(f(x_k+s_k)-\barf(x_k+s_k)\Big) \\*[2ex]
  &      & \hspace*{30mm} + \Big(-f(x_k+s_k)+f(x_k) -\DT_{f,p}(x_k,s_k)\Big)\\*[2ex]
  &      & \hspace*{30mm} + \Big(\DT_{f,p}(x_k,s_k)-\barDT_{f,p}(x_k,s_k)\Big)\Big]\\*[2ex]
  & \leq & \bigfrac{1}{\barDT_{f,p}(x_k,s_k)}
           \left[|f(x_k+s_k)-T_{f,p}(x_k,s_k)|
                             +3\omega |\barDT_{f,p}(x_k,s_k)|\right]\\*[2ex]
  & \leq & \bigfrac{1}{\barDT_{f,p}(x_k,s_k)}
           \left[\bigfrac{L_{f,p}}{(p+1)!}\|s_k\|^{p+1}
                 +{{ 3\omega|\barDT_{f,p}(x_k,s_k)|}}\right]\\*[2ex]
  & \leq & \bigfrac{L_{f,p}}{\sigma_k} + {{3\omega}}\\*[2ex]
  & \leq & 1-\eta.      
  \end{array}
  \]
  Therefore $\rho_k\geq \eta$ and iteration $k$ is successful. 
} 

\subsection{Bounding the expected number of steps with
  $\Sigma_k \geq \sigma_s$}

We now return to the general stochastic process generated by the
\algn algorithm aiming at bounding from above the expected number of steps in
the process generated by the algorithm with $\Sigma_k \geq \sigma_s$. To this
purpose, for all $0 \leq k \leq \ell$, given $\ell \in\iibe{0}{N_{\epsilon}-1}$,
let us define the events 
\begin{eqnarray}
&\Lambda_k=\{\tim{iteration $k$ is such that $\Sigma_k < \sigma_s$}\},
\qquad
\no{\Lambda}_k=\{\tim{iteration $k$ is such that $\Sigma_k \geq \sigma_s$}\}\nonumber\\
&\calS_k=\{\tim{iteration $k$ is successful}\},\nonumber
\end{eqnarray}
and let 
\beqn{defnsigma}
N_\Lambda \eqdef  \sum_{k=0}^{N_\epsilon-1}\indic{\Lambda_k},
\qquad
N_{\no{\Lambda}}
\eqdef \sum_{k=0}^{N_\epsilon-1}\indic{\no{\Lambda}_k},
\eeqn
be the number of steps, in the stochastic process induced by the \algn
algorithm, with $\Sigma_k< \sigma_s$ and $\Sigma_k \geq \sigma_s$, before
iteration $N_{\epsilon}$ is reached, respectively. In what follows we suppose that
AS.1--AS.3 hold.

We may now follow the argument of \cite{CartSche17} to derive an upper bound
on  $\expect\big[N_{\no{\Lambda}}\big]$. In particular, the argument unfolds as follows:
\begin{itemize}
\item[(i)] we apply \cite[Lemma~2.2]{CartSche17} to deduce that, for any
  $\ell\in\iibe{0}{N_\epsilon-1}$ and for all realizations of the \algn
  algorithm, one has that
\beqn{CSlemma22}
\sum_{k=0}^\ell \indic{\no{\Lambda}_k}\indicS \leq \frac{\ell+1}{2};
\eeqn
\item[(ii)] as in \cite{CartSche17}, we note that both
  $\hat\sigma(\indic{\Lambda_k})$ and $\hat\sigma(\indic{\no{\Lambda}_k})$
  belong to $\sigax$, as the random variable $\Lambda_k$ is fully determined
  by the first $k-1$ iterations of the \algn algorithm. Then, setting
  $\ell=N_{\epsilon}-1$  we can rely on \cite[Lemma $2.1$]{CartSche17} (with
$W_k=\indic{\Lambda_k^c}$) and \req{prob-conds} to deduce that  
\beqn{CSlemma21}
\expect\left[ \sum_{k=0}^{N_\epsilon-1}\indic{\no{\Lambda}_k}\indicM \right]
\geq \expect\left[ \sum_{k=0}^{N_\epsilon-1} p_{\calM,k}\indic{\no{\Lambda}_k}\right]
\geq \,\pM\,\expect\left[ \sum_{k=0}^{N_\epsilon-1}  \indic{\no{\Lambda}_k}\right];
\end{equation}
\item[(iii)] as a consequence, given that Lemma \ref{acc-sigma-success-l}
ensures that each iteration $k$ where $\AM_k$ occurs and
$\sigma_k\geq \sigma_s$ is successful, we have that 
\[
\sum_{k=0}^{N_\epsilon-1} \indic{\no{\Lambda}_k}\indicM
\leq  \sum_{k=0}^{N_\epsilon-1} \indic{\no{\Lambda}_k}\indicS
\leq \frac{N_\epsilon}{2},
\]
in which the last inequality follows from \eqref{CSlemma22}, with
$\ell=N_\epsilon-1$. Taking expectation in the above inequality, using
\eqref{CSlemma21} and recalling the rightmost definition in \eqref{defnsigma},
we obtain, as in \cite[Lemma~2.3]{CartSche17}, that, for any realization,
\beqn{bound_Ns}
\expect[N_{\no{\Lambda}}] \leq \frac{1}{2\pM}\expect[N_\epsilon].
\eeqn
\end{itemize} 
\noindent
The remaining upper bound on $\expect[N_\Lambda]$ will be the focus of the next subsection.

\subsection{Bounding the expected number of steps with
  $\Sigma_k < \sigma_s$}

For analyzing $\expect[N_\Lambda]$, where $N_\Lambda$ is defined in
\req{defnsigma}, we now introduce the following variables. 
\vspace{0.1cm}
\begin{definition} With reference to the process \eqref{sprocess} generated by
the \algn algorithm, let us define: 
\vspace{0.1cm}
\begin{equation}
\begin{array}{lcl}
  \bullet\ \overline{\Lambda}_k &\hspace*{-3mm}=&
     \{\mbox{iteration $k$ is such that $\Sigma_k\leq \sigma_s$}\};\\*[1.5ex]
\bullet\ N_I &\hspace*{-3mm}=& \bigsum_{k=0}^{N_{\epsilon}-1}\indic{\overline{\Lambda}_k}\indicMc:
       \mbox{the number of inaccurate iterations with $\Sigma_k \leq \sigma_s$};\\*[1.5ex]
\bullet\ N_A &\hspace*{-3mm}=& \bigsum_{k=0}^{N_{\epsilon}-1}\indic{\overline{\Lambda}_k}\indicM:
       \mbox{the number of accurate iterations with $\Sigma_k \leq \sigma_s$};\\*[1.5ex]
\bullet\ N_{AS} &\hspace*{-3mm}=& \bigsum_{k=0}^{N_{\epsilon}-1}\indic{\overline{\Lambda}_k}\indicM\indicS:
        \mbox{the number of accurate successful iterations with $\Sigma_k \leq \sigma_s$};\\*[1.5ex]
\bullet\ N_{AU} &\hspace*{-3mm}=& \bigsum_{k=0}^{N_{\epsilon}-1}\indic{\Lambda_k}\indicM\indicSc:
        \mbox{the number of accurate unsuccessful iterations with $\Sigma_k < \sigma_s$};\\*[1.5ex]
\bullet\ N_{IS} &\hspace*{-3mm}=& \bigsum_{k=0}^{N_{\epsilon}-1}\indic{\overline{\Lambda}_k}\indicMc\indicS:
        \mbox{the number of inaccurate successful iterations with $\Sigma_k\leq \sigma_s$};\\*[1.5ex]
\bullet\ N_S &\hspace*{-3mm}=& \bigsum_{k=0}^{N_{\epsilon}-1}\indic{\overline{\Lambda}_k}\indicS:
        \mbox{the number of successful iterations with $\Sigma_k\leq \sigma_s$};\\*[1.5ex]
\bullet\ N_U &\hspace*{-3mm}=& \bigsum_{k=0}^{N_{\epsilon}-1}\indic{\Lambda_k}\indicSc:
        \mbox{the number of unsuccessful iterations with $\Sigma_k < \sigma_s$}.
\end{array}
\end{equation}
\vspace{-0.5cm}
\label{def2}
\end{definition}
\noindent
Observe that $\overline{\Lambda}_k$ is the ``closure'' of $\Lambda_k$ in
that the inequality in its definition is no longer strict.

We immediately notice that an upper bound on $\expect[N_\Lambda]$
is available, once an upper bound on $\expect[N_I]+\expect[N_A]$ is
known, since  
\begin{equation}
\label{PlanEN}
\expect[N_\Lambda]
\leq \mathbbm{E}\left[\sum_{k=0}^{N_{\epsilon}-1}\indic{\overline{\Lambda}_k} \right]
= \mathbbm{E}\left[
   \sum_{k=0}^{N_{\epsilon}-1}\indic{\overline{\Lambda}_k}\indicMc
   +\sum_{k=0}^{N_{\epsilon}-1}\indic{\overline{\Lambda}_k}\indicM\right]
= \expect[N_I]+\expect[N_A].
\end{equation}
Using again \cite[Lemma $2.1$]{CartSche17} (with
$W_k=\indic{\overline{\Lambda}_k}$) to give an upper bound on
$\expect[N_I]$, we obtain the following result. 

\llem{a1}{\cite[Lemma~2.6]{CartSche17} 
Let $\AM_k$ be the sequence of events in \eqref{defAMk} and assume that
\req{prob-conds} holds. Let $N_I$, $N_A$ be defined as in Definition
\ref{def2} in the context of the stochastic process \eqref{sprocess}
generated by the \algn algorithm. Then
\begin{equation}
\label{EM1}
\expect[N_I]\le \frac{1-\pM}{\pM} \,\expect[N_A].
\end{equation}
}

\noindent
Turning to  the upper bound for $\expect[N_A]$, we observe that 
\begin{equation}
\label{defEM2}
\expect[N_A]
  =  \expect[N_{AS}]+\expect[N_{AU}]
\leq \expect[N_{AS}]+ \expect[N_U].
\end{equation}
\noindent
Hence, bounding $\expect[N_I]$ can be achieved by providing upper bounds on
$\expect[N_{AS}]$ and $\expect[N_U]$. Regarding the latter, we first
note that the process induced by the \algn algorithm ensures that $\Sigma_k$
is decreased by a factor $\gamma$ on successful steps and increased by the
same factor on unsuccessful ones.
Consequently, by virtue of \cite[Lemma $2,5$]{CartSche17}, we obtain the
following bound. 

\llem{a2}{\cite[Lemma~2.5]{CartSche17}
  For any $\ell\in\{0,...,N_{\epsilon}-1\}$ and for
  all realisations of the \algn algorithm, we have that
\[
\sum_{k=0}^\ell \indic{\Lambda_k}\indicSc
\leq \sum_{k=0}^{\ell}\indic{\overline{\Lambda}_k}\indicS
  +\left\lceil\log_{\gamma}\left( \frac{\sigma_s}{\sigma_0} \right)\right\rceil.
\]
}

\noindent
From this inequality with $\ell=N_{\epsilon}-1$, recalling Definition
\ref{def2} and taking expectations, we therefore obtain that
\begin{equation}
\label{boundEU}
\expect[N_U]
\le \expect[N_S]+ \left\lceil\log_{\gamma}\left( \frac{\sigma_s}{\sigma_0} \right)\right\rceil
=\expect[N_{AS}]+\expect[N_{IS}]
 +\left\lceil\log_{\gamma}\left( \frac{\sigma_s}{\sigma_0} \right)\right\rceil.
\end{equation}
An upper bound on $\expect[N_{AS}]$ is given by the following lemma.

\llem{}{\label{LemmaNAS} Let Assumption AS.1 and AS.2 hold. For all
  realizations of the \algn algorithm we have that 
\begin{equation}
\label{boundEN1}
\expect[N_{AS}]\le \bigfrac{(f_0-\flow)(p+1)!}{{{(\eta-2\omega)}}\sigma_{\min}\psi(\sigma_s)}
            \left(\min_{j\in\ii{q}}\epsilon_j\right)^{-\varpi}+1,
\end{equation}
where $\varpi$, $\psi(\sigma)$ and $\sigma_s$ are defined in \req{pi-def},
\req{psi-def} and \req{large-sigma}, respectively.
}

\proof{
For all realizations of the \algn algorithm we have that:
\begin{itemize}
\item if iteration $k$ is successful, then \eqref{acc-f-decrease-1} holds;
\item if iteration $k$ is successful and accurate (i.e.,
  $\indicS\indicM=1$) and \req{strong} fails
  for $(x_{k+1},\delta_{k,j})$, then \req{longs} holds;
\item if iteration $k$ is unsuccessful, the mechanism of the \algn algorithm
  guarantees that $x_k=x_{k+1}$ and, hence, that $f(x_{k+1})=f(x_k)$.
\end{itemize}
Therefore, for any $\ell\in\{0,...,N_{\epsilon}-1\}$,
\begin{eqnarray}
f_0-\flow&\ge& f_0-f(X_{\ell+1})
=\sum_{k=0}^{\ell}\indicS(f(X_k)-f(X_{k+1}))
\geq \sum_{k=0}^{\ell}\indicS \frac{{{(\eta-2\omega)}}\sigma_{\min}}{(p+1)!}\|S_k\|^{p+1}\nonumber\\
& \geq & \sum_{k=0}^{\ell-1}\indicS\indicM
       \frac{{{(\eta-2\omega)}}\sigma_{\min}}{(p+1)!}\|S_k\|^{p+1}\\
& \geq & \sum_{k=0}^{\ell-1} \indicS\indicM
      \frac{{{(\eta-2\omega)}}\sigma_{\min}}{(p+1)!} \psi(\Sigma_k)
      \left(\min_{j\in\ii{q}}\epsilon_j\right)^\varpi\nonumber\\
& \geq & \sum_{k=0}^{\ell-1}
      \indicS\indicM\indic{\overline{\Lambda}_k}
      \frac{{{(\eta-2\omega)}}\sigma_{\min}}{(p+1)!} \psi(\Sigma_k)
      \left(\min_{j\in\ii{q}}\epsilon_j\right)^\varpi\nonumber\\
& \geq & \frac{{{(\eta-2\omega)}}\sigma_{\min}}{(p+1)!} \psi(\sigma_s)
      \left(\min_{j\in\ii{q}}\epsilon_j\right)^\varpi
      \left(\sum_{k=0}^{\ell-1} \indicS \indicM
      \indic{\overline{\Lambda}_k}\right),\label{last}
\end{eqnarray}
having set $f_0\eqdef f(X_0)$ and where the last inequality is due to fact
that $\psi(\sigma)$ is a decreasing function. We now notice that, by Definition \ref{def2},
\[
N_{AS}-1\le \sum_{k=0}^{N_{\epsilon}-2} \indic{\overline{\Lambda}_k} \indicM \indicS.
\]
Hence, letting $\ell=N_{\epsilon}-1$ and taking expectations in \eqref{last}, we conclude that
\[
f_0-\flow \ge (\expect[N_{AS}]-1)\frac{{{(\eta-2\omega)}}\sigma_{\min}}{(p+1)!}
\psi(\sigma_s) \left(\min_{j\in\ii{q}}\epsilon_j\right)^\varpi,
\]
which is equivalent to \eqref{boundEN1}.
}
\noindent
While inequalities \req{boundEN1} and \req{boundEU}   provide upper bounds on $\expect[N_{AS}]$
and $\expect[N_U]$, as desired, the latter still 
depends on $\expect[N_{IS}]$, which has to be bounded from above as well. This
can be done by following \cite{CartSche17} once more: 
Definition \ref{def2}, \eqref{EM1} and \eqref{defEM2} directly imply that 
\begin{equation}
\label{1boundEM3}
\expect[N_{IS}]
\leq \expect[N_I]
\leq  \frac{1-\pM}{\pM} \expect[N_A]
\leq \frac{1-\pM}{\pM}\left(\expect[N_{AS}] +  \expect[N_U] \right)
\end{equation}
and hence
\begin{equation}
\label{boundEM3}
\expect[N_{IS}]
\leq \frac{1-\pM}{2\pM-1}\left(2\expect[N_{AS}]
      + \left\lceil\log_{\gamma}\left( \frac{\sigma_s}{\sigma_0} \right)\right\rceil\right)
\end{equation}
follows from \eqref{boundEU} (remember that $\half < \pM \leq 1$). Thus, the right-hand side in \eqref{EM1} is in
turn upper bounded by virtue of \eqref{defEM2}, \eqref{boundEU},
\eqref{boundEM3} and \eqref{boundEN1}, giving
\begin{eqnarray}
\expect[N_A]
& \leq & \expect[N_{AS}]+\expect[N_U]
\leq 2\expect[N_{AS}]+\expect[N_{IS}]
    +\left\lceil\log_{\gamma}\left(\frac{\sigma_s}{\sigma_0}\right)\right\rceil\nonumber\\
& \leq & \left(\frac{1-\pM}{2\pM-1}+1\right)\left(2\expect[N_{AS}]
      + \left\lceil\log_{\gamma}\left( \frac{\sigma_s}{\sigma_0} \right)\right\rceil\right)\nonumber\\
&   =  & \frac{\pM}{2\pM-1}\left[2\expect[N_{AS}]
      + \left\lceil\log_{\gamma}\left( \frac{\sigma_s}{\sigma_0} \right)\right\rceil\right]\nonumber\\
& \leq & \frac{\pM}{2\pM-1}\left[\bigfrac{2(f_0-\flow)(p+1)!}{{{(\eta-2\omega)}}\sigma_{\min}\psi(\sigma_s)}
        \left(\min_{j\in\ii{q}}\epsilon_j\right)^{-\varpi}
        + \left\lceil\log_{\gamma}\left( \frac{\sigma_s}{\sigma_0} \right)\right\rceil+2\right]. \label{EM2}
\end{eqnarray}
This inequality, together with \eqref{PlanEN} and \eqref{EM1}, finally
gives the desired bound on $\expect[N_{\Lambda}]$:
\begin{equation}
\label{bound_Nsc}
\expect[N_{\Lambda}]
\leq \frac{1}{\pM}\expect[N_A]
\leq \frac{1}{2\pM-1}\left[\bigfrac{2(f_0-\flow)(p+1)!}{{{(\eta-2\omega)}}\sigma_{\min}\psi(\sigma_s)}
  \left(\min_{j\in\ii{q}}\epsilon_j\right)^{-\varpi}
  + \left\lceil\log_{\gamma}\left( \frac{\sigma_s}{\sigma_0} \right)\right\rceil+2\right].
\end{equation}

\noindent
We can now express our final complexity result in full.

\lthm{complexity-th}{
Suppose that AS.1--AS.3 hold. Then the following
conclusions also hold.
\begin{enumerate}
\item If $q\in\{1,2\}$,  then
\[
\expect[N_\epsilon]
\leq \kappa(\pM)
\left(\frac{2(f_0-\flow)(p+1)!}{{{(\eta-2\omega)}}\sigma_{\min}\psi(\sigma_{s})}
  \left(\min_{j\in\ii{q}}\epsilon_j\right)^{-\frac{p+1}{p-q+1}}
  +\left\lceil\log_{\gamma}\left( \frac{\sigma_s}{\sigma_0}\right)\right\rceil+2 \right),
\]
\item If $q > 2$, then
\[
\expect[N_\epsilon]
\leq \kappa(\pM)
\left(\frac{2(f_0-\flow)(p+1)!}{{{(\eta-2\omega)}}\sigma_{\min}\psi(\sigma_{s})}
 \left(\min_{j\in\ii{q}}\epsilon_j\right)^{-\frac{q(p+1)}{p}}
  +\left\lceil\log_{\gamma}\left( \frac{\sigma_s}{\sigma_0}\right)\right\rceil +2\right),
\]
\end{enumerate}
with $\kappa(\pM) \eqdef \frac{2\pM}{(2\pM-1)^2}$ and $N_\epsilon$, $\psi(\sigma)$, $\sigma_s$ defined as in
\req{Neps-def}, \req{psi-def}, \req{large-sigma}, respectively.
}

\noindent
\proof{
Recalling the definitions \eqref{defnsigma} and the bound \eqref{bound_Ns}, we obtain that
\[
\expect[N_{\epsilon}]=\expect[N_{\Lambda}^c]+\expect[N_{\Lambda}]
\leq \frac{\expect[N_{\epsilon}]}{2\pM}+\expect[N_{\Lambda}],
\]
which implies, using \req{bound_Nsc}, that
\[
\frac{2\pM-1}{2\pM}\expect[N_{\epsilon}]
\leq\frac{1}{2\pM-1}\left(\bigfrac{2(f_0-\flow)(p+1)!}{{{(\eta-2\omega)}}\sigma_{\min}\psi(\sigma_s)}
\left(\min_{j\in\ii{q}}\epsilon_j\right)^{-\varpi}
+ \left\lceil\log_{\gamma}\left( \frac{\sigma_s}{\sigma_0} \right)\right\rceil+2\right).
\]
This bound and the inequality $\half < \pM \leq 1$ yield the desired result.
}

\noindent
Since the \algn algorithm requires at most two function evaluations and one
evaluation of the derivatives of orders one to $p$ per 
iteration, the bounds stated in the above theorem effectively provide an upper
bound on the average evaluation complexity of finding
$(\epsilon,\delta)$-approximate $q$-th order
minimizers.

Theorem~\ref{complexity-th} generalizes the complexity bounds
stated in \cite[Theorem~5.5]{CartGoulToin20a} to the case where evaluations of
$f$ and its derivatives are inexact, under probabilistic assumptions on the
accuracies of the latter. Remarkably, the bounds of
Theorem~\ref{complexity-th} \emph{ are essentially identical in order of the tolerance $\epsilon$
to those obtained in \cite[Theorem~5.5]{CartGoulToin20a},} in that they only
differ by the presence of an additional term in
$|\log(\min_{j\in\ii{q}}\epsilon_j)|$.  Moreover, it was shown in
\cite[Theorems~6.1 and 6.4]{CartGoulToin20a} that the evaluation complexity
bounds are sharp in $\epsilon$ for exact evaluations
and Lipschitz continuous derivatives of $f$.  Since the \algn algorithm
reduces to the algorithm proposed in that reference when all values are exact,
we deduce that the lower bound on evalution complexity presented in this
reference is also valid in our case.  Thus, considering that, for small
$\epsilon_j$, the term $|\log(\min_{j\in\ii{q}}\epsilon_j)|$ is much smaller that the terms in
$\min_{j\in\ii{q}}\epsilon_j^{-(p+1)/(p-q+1)}$ or
$\min_{j\in\ii{q}}\epsilon_j^{-q(p+1)/p}$,
we conclude that the \emph{the presence of
random noise in the derivatives and of inexactness in function values does not
affect the evaluation complexity} of adaptive regularization algorithms for
the local solution of problem \req{problem}.  In addition, we also deduce that
\emph{the complexity bounds of Theorem~\ref{complexity-th} are
essentially\footnote{Modulo the negligible logarithmic term.} sharp in order of $\epsilon$.}

It is interesting to compare our results with those of
\cite{Arjeetal20}. These authors mention an ``elbow effect'' for algorithms
using randomly perturbed derivatives in that they state a lower bound on
evaluation complexity for second-order approximate minimizers of
$\mathcal{O}(\epsilon_2^{-3})$ for all $p \ge 2$, in contrast with our
smoothly decreasing $\mathcal{O}(\min[\epsilon_1,\epsilon_2]^{-(p+1)/(p-1)})$
bound.  However, their framework is very different. Firstly, they assume the a
priori knowledge of the Lipschitz constants, which makes monitoring of the
function values unnecessary in an adaptive regularization algorithm, an
assumption we have explictly avoided for consistency. Most importantly, their
accuracy model is significantly more permissive than ours, as it
allows\footnote{See definitions (76) and (87) in \cite{Arjeetal20}.}
derivatives' estimates of the form $\overline{\nabla_x^j f}(x) = (z/\mu)
\nabla_x^j f(x)$ where $z$ is a $(0,1)$ Bernoulli random variable of
parameter $\mu$ . Although unbiased and of bounded variance under AS.1,
such estimates result, for nonzero $\nabla_x^j f(x)$, in an infinite relative
error with probability $1-\mu$.  Since they consider values of $\mu$ of the
order of  $\epsilon_2^2$, this is clearly too loose for AS.3 to hold. This
illustrate that, unsurprisingly, the evaluation complexity bound for
algorithms using inexact information strongly depends on the specific
(potentially probabilistic) accuracy model considered.

We conclude this section by noting that the complexity bounds we have derived
depend on the smallest of the accuracy thresholds $\epsilon_j$.  We could
therefore derive the complete theory with a single $\epsilon$ for all
optimality orders, marginally improving notation.  We have refrained from
doing so because users of numerical optimization algorithms very rarely
makes this choice in practice, but typically uses application- and
order-dependent thresholds.

\numsection{Extension to convexly constrained problems}

As it turns out, it is easy to extend the above results to the case where the
problem is convexly constrained, that is when \req{problem} is replaced by
\beqn{problem-convex}
\min_{x \in \calX}  f(x), 
\eeqn
where $\calX$ is a convex subset of $\Re^n$. We have refrained from
considering this problem from the start for clarity of exposition, but we now
review the (limited) changes that are necessary to cover this more general
problem. 
\begin{enumerate}
\item We may first weaken AS.1 to require that $f$ is $p$ times continuously
  differentiable in an open convex neighbourhood of $\calX$ and that the
  Lipschitz conditions \req{f-holder} only hold in that neighbourhood.
\item We must then revise our approximate criticality measure \req{phi-def}
  to reflect the constrained nature of \req{problem-convex}. This is done
  by considering the Taylor decrement which is achievable only for
  displacements $d$ which preserve feasibility. We may therefore replace
  $\phi_{f,j}^{\delta_j}(x)$ in \req{phi-def} by
  \[
  \phi_{f,j}^{\delta_j}(x) = \max_{x+d \in \calX, \,\|d\| \leq \delta_j} \DT_{f,j}(x,d)
  \]
  for $x\in \calX$. This new definition is then used in \req{Neps-def} to
  obtain a new family of stopping times.
\item We next insist that feasibility is maintained throughout the execution
  of the algorithm, in that we require that $x_0\in\calX$ and that $s_k$ is
  computed such that the trial point $x_k+s_k$ is also feasible. Moreover, our
  criterion for terminating the step search must also reflect its constrained
  nature, which is obtained by replacing $\barphi_{m_k,j}^{\delta_{k,j}}(s_k)$
  in \req{step-term} by
  \[
  \barphi_{m_k,j}^{\delta_{k,j}}(s_k)
  = \max_{x_k+s_k+d\in\calX, \,|d\|\leq \delta_{k,j}}\barDT_{m_k,j}(s_k,d)
  \]
  for $x_k+s_k\in \calX$.
\item The theory is then unchanged for this new context, with one
  caveat. We note that the proof of Lemma~\ref{delta-lower-l} for the case
  where $q\in\{1,2\}$ does depend on the fact that $T_{m_k,j}(s_k^*,d)$ is a
  convex function of $d$ for $j\in\{1,2\}$ because of the unconstrained
  optimality conditions. Obviously, while maintaining convexity is possible in the
  convexly constrained case when $q=1$, it may now fail for $q=2$. As a
  consequence, this case must be considered in the same way as for other
  larger values of $q$.  This then imposes that we have to change the
  condition ``if $q\in\{1,2\}$''  or ``if $j\in\{1,2\}$'' to ``if $q=1$'' or
  ``if $j=1$'', respectively,   in \req{defAMk}, the first statement of
  Lemma~\ref{delta-lower-l}, \req{pi-def} and its proof, and in the first
  statement of Theorem~\ref{complexity-th}.
\end{enumerate}
It is remarkable that no further change is necessary for deducing
Theorem~\ref{complexity-th} for problem \req{problem-convex}.
This extension to the convexly convex case is also a novel feature for
algorithms considering randomly perturbed derivatives.

\numsection{Conclusions and perspectives}\label{concl-s}

We have shown that the \algn algorithm, a stochastic inexact adaptive
regularization algorithm using derivatives of order up to $p$, computes an 
$(\epsilon,\delta)$-approximate $q$-th order minimizer of $f$ in
problem~\req{problem} in at most $O(\epsilon^{-\frac{p+1}{p-q+1}})$
iterations in expectation if $q$ is either one or two,
while it may need $O(\epsilon^{-\frac{q(p+1)}{p}})$ iterations in
expectation in the other cases\footnote{These simplified order bounds assume that
$\epsilon_j=\epsilon$ for $j\in\ii{q}$.}. Moreover, these bounds are
essentially sharp in
the order of $\epsilon$ (see \cite{CartGoulToin20a}). We therefore conclude
that, if the probabilities $p_{\calM,k}$ in AS.3 are suitably large, the evaluation
complexity of the \algn algorithm is identical (in order) to that of the exact
algorithm in \cite{CartGoulToin20a}. We finally provided an extension of these
results to the convexly constrained case.

We also note that the full power of AS.1 is only required for
Lemma~\ref{delta-lower-l}, while Lipschitz continuity of $\nabla_x^p f(x)$ is
sufficient for all subsequent derivations.  Thus if suitable lower bounds on
$\Delta_{k,j}$ can be ensured in some other way, our development remains valid
(although the precise complexity bounds will depend on the new bounds on
$\Delta_{k,j}$). In AS.1, we have also required (Lipschitz) continuity of $f$
and its derivatives in $\Re^n$. This can be weakened to  requiring this
property only on the ``tree of iterates'' $\cup_{k\geq0}[x_k,x_k+s_k]$, but
this weaker assumption is often impossible to verify \emph{a priori}.
In the same vein,
it also is possible to avoid requiring that \req{delta-lower} is always
ensured by the \algn algorithm whenever $q>2$ by instead redefining $\AM_k$ to
also include the satisfaction of this condition. We have preferred using an
explicit assumption because this approach better differentiates deterministic
requirements on the algorithm from stochastic assumptions more related to the
problem itself.

We finally recall that \cite{CartGoulToin20a} also derives complexity bounds
for the (possibly non-smooth) composite optimization problem.  We expect that
the theory presented here can be extended to also cover this case.

An analysis covering adaptive regularization algorithms where the objective
function evaluations are also subject to general random noise,
parallel to that provided for trust-region methods for low order minimizers
in \cite{BlanCartMeniSche19}, remains, for now, an open and challenging
question.

{\footnotesize
\section*{\footnotesize Acknowledgment}

INdAM-GNCS partially supported the first and third authors under Progetti di
Ricerca 2019. The second author was partially supported by INdAM through a
GNCS grant. The last author gratefully acknowledges the support and
friendly environment provided by the Department of Industrial Engineering at
the Universit\`{a} degli Studi di Firenze (Italy) during his visit in the fall
of 2019.


\section*{Appendix}

{\bf Proof of Lemma~\ref{step-exists}}
\setcounter{equation}{0}
\renewcommand{\theequation}{A.\arabic{equation}}
\renewcommand{\thesection}{A}
    
Let $s_k^*$ be a global minimizer of $m_k(s)$. By Taylor's theorem, we have that, for all
$d$,
\beqn{1star}
\begin{array}{lcl}
0 & \leq &  m_k(s_k^*+d)-m_k(s_k^*)
= \bigsum_{\ell=1}^p \bigfrac{1}{\ell!}\nabla_s^\ell \barT_{f,p}(x_k,s_k^*)[d]^\ell\\*[2.5ex]
& & \hspace*{5mm} + \bigfrac{\sigma_k}{(p+1)!}
\Bigg[\bigsum_{\ell=1}^p\frac{1}{\ell !}\nabla_s^\ell \left(\|s_k^*\|^{p+1}\right)[d]^\ell
       + \frac{1}{(p+1)!}\nabla_s^{p+1} \left(\|s_k^*+\tau d\|^{p+1}\right)[d]^{p+1} \Bigg]
\end{array}
\eeqn
for some $\tau \in (0,1)$.
We may now use the expression of $\nabla_s^\ell \left(\|s_k^*\|^{p+1}\right)$
given by \cite[Lemma~2.4]{CartGoulToin20b} in \req{1star} 
and deduce that, for any $j\in \ii{q}$ and all $d$,
\beqn{appext-comp-sl-1}
\begin{array}{l}
- \bigsum_{\ell=1}^j \bigfrac{1}{\ell!}\nabla_s^\ell \barT_{f,p}(x_k,s_k^*)[d]^\ell
- \bigfrac{\sigma_k}{(p+1)!}\bigsum_{\ell=1}^j\nabla_s^\ell \|s_k^*\|^{p+1}[d]^\ell\\*[2.5ex]
\hspace*{15mm}
\leq \bigsum_{\ell=j+1}^p \bigfrac{1}{\ell!}\nabla_s^\ell \barT_{f,p}(x_k,s_k^*)[d]^\ell
     + \bigfrac{\sigma_k}{(p+1)!}
     \Bigg[\bigsum_{\ell=j+1}^p\frac{1}{\ell !}\nabla_s^\ell \|s_k^*\|^{p+1}[d]^\ell
            + \|d\|^{p+1}\Bigg].
\end{array}
\eeqn
It is now possible to choose $\delta_{k,j}\in
(0,1]$ such that, for every $d$ with $\|d\| \leq \delta_{k,j}$,
\beqn{appext-comp-sl-2}
\begin{array}{l}  
\bigsum_{\ell=j+1}^p \bigfrac{1}{\ell!}\nabla_s^\ell \barT_{f,p}(x_k,s_k^*)[d]^\ell
      + \bigfrac{\sigma_k}{(p+1)!}
        \Bigg[\bigsum_{\ell=j+1}^p\bigfrac{1}{\ell!}\nabla_s^\ell \|s_k^*\|^{p+1}[d]^\ell
               + \|d\|^{p+1} \Bigg]\\*[2ex]
\hspace*{30mm} \leq \bigfrac{1}{2}\theta \epsilon_j \,\bigfrac{\delta_{k,j}^j}{j!}.
\end{array}
\eeqn
We therefore obtain that if $\delta_{k,j}$ is small enough to ensure
\req{appext-comp-sl-2}, then \req{appext-comp-sl-1} implies that
\beqn{roundstar}
- \bigsum_{\ell=1}^j \bigfrac{1}{\ell!}\nabla_s^\ell \barT_{f,p}(x_k,s_k^*)[d]^\ell
- \bigfrac{\sigma_k}{(p+1)!}\bigsum_{\ell=1}^j\nabla_s^\ell \|s_k^*\|^{p+1}[d]^\ell
\leq \frac{1}{2} \theta \epsilon_j\,\bigfrac{\delta_{k,j}^j}{j!}.
\eeqn
and therefore that, for all $j\in\ii{q}$,
\[
\max_{\|d\|\leq \delta_{k,j}}\barDT_{m_k,j}(s_k^*,d)
\leq \frac{1}{2}\theta \epsilon_j\,\bigfrac{\delta_{k,j}^j}{j!}.
\]
Thus the pair $(s_k^*,\delta_k)$ is acceptable for Step~2 of the algorithm.
If we assume now that $x_k+s_k^*$ is not an isolated feasible point,
the above inequality and continuity of
$\barT_{f,p}(x_k,s)$ and its derivatives with respect to $s$ then ensure the
existence of a feasible neighbourhood $\calN_k^*$ of $s_k^*$ in which 
\beqn{unc-term}
\max_{\|d\|\leq \delta_{k,j}}\barDT_{m_k,j}(s,d)
\leq \theta \epsilon_j\,\bigfrac{\delta_{k,j}^j}{j!}.
\eeqn
for all $s \in \calN_k^*$. 
We may then choose any $s_k$ in $\calN_k^*$ such that, in addition to
satisfying \req{unc-term} and being such that $x_k+s_k$ is
feasible, \req{descent} also holds.  Thus the definition of
$\barphi_{m_k,j}^{\delta_{k,j}}(s_k)$ in \req{step-term} gives that
\beqn{part1phim}
\barphi_{m_k,j}^{\delta_{k,j}}(s_k)
\leq \theta \epsilon_j\,\bigfrac{\delta_{k,j}^j}{j!}
\eeqn
and every such $(s_k,\delta_k)$ is also acceptable for Step~2 of the algorithm.
\vskip 10 pt
{\bf Proof of Lemma~\ref{delta-lower-l}}

Let $s_k^*$ be a global minimizer of $m_k(s)$.   
We first consider the case where $q\in\{1,2\}$.  Then it
is easy to verify that, for each $j\in \ii{q}$, the optimization problem
involved in the definition of $\barphi_{m_k,j}^{\delta_{k,j}}(s_k^*)$ (in
\req{step-term}) is convex and therefore that $\delta_{k,j}$ can be chosen
arbitrarily in $(0,1]$.  The first case of Lemma~\ref{delta-lower-l} then
follows from the continuity of $\barphi_{m_k,j}^{\delta_{k,j}}(s)$ with
respect to $s$. Unfortunately, the crucial convexity property is lost for $q>2$

Unfortunately, the crucial convexity property is lost for $q>2$ and, in order
to prove the second case, we now pursue the reasoning of the proof of
Lemma~\ref{step-exists}. We
start by supposing that $\|s_k^*\| > 1$. We may then reduce the neighbourhood of
$s_k^*$ in which $s_k$ can be chosen enough to guarantee that
$\|s_k\| \geq 1$,
which then gives the desired result because of \req{unc-term}.
Suppose therefore that $\|s_k^*\| \leq 1$. The triangle
inequality then implies that \[
\|\nabla_s^\ell \barT_{f,p}(x_k,s_k^*)\|
\leq \bigsum_{i=\ell}^p\bigfrac{1}{(i-\ell)!}\|\overline{\nabla_x^i f}(x_k)\|\,\|s_k^*\|^{i-\ell},
\]
for $\ell \in \iibe{q+1}{p}$, and thus,
using, AS.1 and \cite[Lemma~2.4]{CartGoulToin20b}, we deduce that
\[
\begin{array}{l}
\bigsum_{\ell=j+1}^p \bigfrac{1}{\ell!}\nabla_s^\ell \barT_{f,p}(x_k,s_k^*)[d]^\ell
     + \bigfrac{\sigma_k}{(p+1)!}
       \Bigg[\bigsum_{\ell=j+1}^p\nabla_s^\ell \|s_k^*\|^{p+1}[d]^\ell \Bigg]\\*[3ex]
\hspace*{20mm} \leq  \bigsum_{\ell=j+1}^p\bigfrac{\|d\|^\ell}{\ell!}
       \Bigg[\bigsum_{i=\ell}^p\bigfrac{\|s_k^*\|^{i-\ell}}{(i-\ell)!}
       \|\overline{\nabla_x^i f}(x_k)\|+\bigfrac{\sigma_k\|s_k^*\|^{p-\ell+1}}{(p-\ell+1)!}\Bigg].
\end{array}
\]
We now call upon the fact that, since $q \geq 3$ and $\AM_k$ occurs by
assumption, $\AM_k^{(4)}$ also occurs.  Thus
\[
\begin{array}{l}
\bigsum_{\ell=j+1}^p \bigfrac{1}{\ell!}\nabla_s^\ell \barT_{f,p}(x_k,s_k^*)[d]^\ell
     + \bigfrac{\sigma_k}{(p+1)!}
       \Bigg[\bigsum_{\ell=j+1}^p\nabla_s^\ell \|s_k^*\|^{p+1}[d]^\ell \Bigg]\\*[3ex]
       \hspace*{20mm} \leq  \bigsum_{\ell=j+1}^p \bigfrac{\|d\|^\ell}{\ell!}
       \Bigg[ \Theta \, \bigsum_{i=\ell}^p\bigfrac{\|s_k^*\|^{i-\ell}}{(i-\ell)!}
         + \bigfrac{\sigma_k\|s_k^*\|^{p-\ell+1}}{(p-\ell+1)!}\Bigg].
\end{array}
\]
We therefore obtain from \req{appext-comp-sl-2} that any pair $(s_k^*,\delta_{s,j})$
satisfies \req{roundstar}  for $\|d\|\leq \delta_{s,j}$ if
\beqn{unc3-arqp-accept-sstar}
\bigsum_{\ell=j+1}^p \bigfrac{\delta_{s,j}^\ell}{\ell!} \left[\Theta \bigsum_{i=\ell}^p
  \bigfrac{1}{(i-\ell)!}\|s_k^*\|^{i-\ell}
  + \bigfrac{\sigma_k\|s_k^*\|^{p-\ell+1}}{(p-\ell+1)!}\right]
  + \sigma_k \bigfrac{\delta_{s,j}^{p+1}}{(p+1)!}
\leq \frac{1}{2} \theta\epsilon_j\,\bigfrac{\delta_{s,j}^j}{j!}.
\eeqn
which, because $\|s_k^*\|\leq 1$, is in turn ensured by the inequality
\beqn{unc3-arqp-conds1}
\sum_{\ell=j+1}^p \frac{\delta_{s,j}^\ell}{\ell!}
     \left[ \Theta\,\sum_{i=\ell}^p\frac{1}{(i-\ell)!}+\sigma_k\right]
     + \sigma_k \bigfrac{\delta_{s,j}^{p+1}}{(p+1)!}
\leq \frac{1}{2}\theta \epsilon_j\,\frac{\delta_{s,j}^j}{j!}.
\eeqn
Observe now that, since $\delta_{s,j}\in [0,1]$, $\delta_{s,j}^\ell \leq
\delta_{s,j}^{j+1}$ for $\ell\in\iibe{j+1}{p}$.  Moreover, we have that,
\[
\sum_{i=\ell}^p\frac{1}{(i-\ell)!} \leq e < 3, \ms(\ell \in \iibe{j+1}{p+1}),
\ms\ms
\sum_{\ell=j+1}^{p+1}\frac{1}{\ell!} \leq e-1 < 2
\]
and therefore \req{unc3-arqp-conds1} is guaranteed by the condition
\beqn{unc3-arqp-safecond-L}
j!(6\Theta+2\sigma_k) \,\delta_{s,j}
\leq  \frac{1}{2}\theta\epsilon_j,
\eeqn
which means that the pair $(s_k^*,\delta_s)$ satisfies \req{roundstar} for
all $j \in \ii{q}$ whenever,
\[
\delta_{s,j}
\leq \frac{1}{2} \delta_{\min,k}
\eqdef \frac{\theta\epsilon_j}
                     {2q!(6\Theta+2\sigma_k)}.
                     \]
As in the proof of Lemma~\ref{step-exists}, we may invoke continuity of the
derivatives of $m_k(s)$ with respect to $s$ to
deduce that there exists a neighbourhood $\calN_k^*$ of $s_k^*$  such that
\req{unc-term}  holds for every $s\in \calN_k^*$ and every $\delta_{k,j}\leq \delta_{\min,k}$.
Choosing now $s_k$ to ensure \req{descent} in addition
to \req{unc-term}, we obtain that
the pair $(s_k,\delta_{k,j})$ satisfies both \req{descent} and
\[
\barphi_{m_k,j}^{\delta_{k,j}}(s_k)
\leq\theta \epsilon_j\,\frac{\delta_{k,j}^j}{j!}.
\]
The desired conclusion then follows with
\[
\kappa_\delta(\sigma)
= \frac{\nu\theta}{q!(6\Theta+2\sigma)}
\]
for any constant $\nu \in (0,1)$. Moreover, $\kappa_\delta(\sigma)$ is clearly
a decreasing function of $\sigma$.

\end{document}